\RequirePackage{etoolbox}\csdef{input@path}{{style/}{graphics/}}
\documentclass[aos,MSNbibl,nameyear,dvips]{arximspdf}
\usepackage{multirow,dcolumn}
\usepackage{float}

\doi{10.1214/13-AOS1198} 
\volume{42}
\issue{3}
\pubyear{2014}
\firstpage{819}
\lastpage{849}

\makeatletter
\newcolumntype{d}[1]{D{.}{.}{#1}}
\newproclaim{rem}{Remark}
\newproclaim{mo}{Model}
\newtheorem{theorem}{Theorem}
\newtheorem{lemma}{Lemma}
\newtheorem{corollary}{Corollary}
\floatstyle{ruled}
\newfloat{algorithm}{tbhp}{loa}
\floatname{algorithm}{Algorithm}
\makeatother

\begin{document}
\begin{frontmatter}

\title{Strong oracle optimality of folded concave penalized estimation}
\runtitle{On the computable strong oracle optimality}

\begin{aug}
\author[A]{\fnms{Jianqing} \snm{Fan}\ead[label=e1]{jqfan@princeton.edu}\thanksref{t1}},
\author[B]{\fnms{Lingzhou} \snm{Xue}\ead[label=e2]{lzxue@psu.edu}\thanksref{t2}}
\and
\author[C]{\fnms{Hui} \snm{Zou}\corref{}\ead[label=e3]{zouxx019@umn.edu}\thanksref{t3}}
\runauthor{J. Fan, L. Xue and H. Zou}
\affiliation{Princeton University, Pennsylvania State University and\break University of Minnesota}
\address[A]{J. Fan\\
Department of Operations Research\\
\quad and Financial Engineering\\
Princeton University\\
Princeton, New Jersey 08544\\
USA\\
\printead{e1}}
\address[B]{L. Xue\\
Department of Statistics\\
Pennsylvania State University\\
University Park, Pennsylvania 16802\\
USA\\
\printead{e2}}
\address[C]{H. Zou\\
School of Statistics\\
University of Minnesota\\
Minneapolis, Minnesota 55414\\
USA\\
\printead{e3}}
\end{aug}
\thankstext{t1}{Supported by NIH Grant R01-GM072611 and NSF Grants
DMS-12-06464 and 0704337.}
\thankstext{t2}{Supported by NIH Grant R01-GM100474 as a postdoctor at
Princeton University.}
\thankstext{t3}{Supported by NSF Grant DMS-08-46068 and a grant from
Office of Naval Research.}

\received{\smonth{10} \syear{2012}}
\revised{\smonth{12} \syear{2013}}

%
\begin{abstract}
Folded concave penalization methods have been shown to enjoy the strong
oracle property for high-dimensional sparse estimation. However,
a~folded concave penalization problem usually has multiple local
solutions and the oracle property is established only for one of the
unknown local solutions. A challenging fundamental issue still remains
that it is not clear whether the local optimum computed by a given
optimization algorithm possesses those nice theoretical properties. To
close this important theoretical gap in over a decade, we provide a
unified theory to show explicitly how to obtain the oracle solution via
the local linear approximation algorithm. For a folded concave
penalized estimation problem, we show that as long as the problem is
localizable and the oracle estimator is well behaved, we can obtain the
oracle estimator by using the one-step local linear approximation. In
addition, once the oracle estimator is obtained, the local linear
approximation algorithm converges, namely it produces the same
estimator in the next iteration. The general theory is demonstrated by
using four classical sparse estimation problems, that is, sparse linear
regression, sparse logistic regression, sparse precision matrix
estimation and sparse quantile regression.
\end{abstract}

%

\end{frontmatter}

\section{Introduction}\label{sec1}

Sparse estimation is at the center of the stage of high-dimensional
statistical learning. The two mainstream methods are the LASSO (or~$\ell_1$ penalization) and the folded concave penalization
[\citet{scad}] such as the SCAD and the MCP. Numerous papers have
been devoted to the numerical and theoretical study of both methods. A
strong irrepresentable condition is necessary for the LASSO to be
selection consistent [\citet{meinshausen2006}, \citet
{yu2006}, \citet{alasso}]. The folded concave penalization,
unlike the LASSO, does not require the irrepresentable condition to
achieve the variable selection consistency and can correct the
intrinsic estimation bias of the LASSO [\citet{scad}, \citet
{fan2004}, \citet{mcp}, \citet{fan2011}]. The LASSO owns
its popularity largely to its computational properties. For certain
learning problems, such as the LASSO penalized least squares, the
solution paths are piecewise linear which allows one to employ a
LARS-type algorithm to compute the entire solution path efficiently
[\citet{lars}]. For a more general class of $\ell_1$
penalization problems, the coordinate descent algorithm has been shown
to be very useful and efficient [\citeauthor{friedman2008} (\citeyear{friedman2008,friedman2010})].

The computation for folded concave penalized methods is much more
involved, because the resulting optimization problem is usually
nonconvex and has multiple local minimizers. Several algorithms have
been developed for computing the folded concave penalized estimators.
\citet{scad} worked out the local quadratic approximation (LQA)
algorithm as a unified method for computing the folded concave
penalized maximum likelihood. \citet{onestep} proposed the local
linear approximation (LLA) algorithm which turns a concave penalized
problem into a series of reweighed $\ell_1$ penalization problems.
Both LQA and LLA are related to the MM principle [\citet
{hunter2004}, \citet{hunter2005}]. Recently, coordinate descent
was applied to solve the folded concave penalized least squares
[\citet{sparsenet}, \citet{fan2011}]. \citet{mcp}
devised a PLUS algorithm for solving the penalized least squares using
the MCP and proved the oracle property. \citeauthor{tzhang2010}
(\citeyear{tzhang2010,tzhang2013}) analyzed the capped-$\ell_1$
penalty for solving the penalized least squares and proved the oracle
property as well. With these advances in computing algorithms, one can
now at least efficiently compute a local solution of the folded concave
penalized problem. It has been shown repeatedly that the folded concave
penalty performs better than the LASSO in various high-dimensional
sparse estimation problems. Examples include sparse linear regression
[\citet{scad}, \citet{mcp}], sparse generalized linear
model [\citet{fan2011}], sparse Cox's proportional hazards model
[\citet{fan2012}], sparse precision matrix estimation
[\citet{fan2009a}], sparse Ising model [\citet
{xue-zou-cai2012}], and sparse quantile regression [\citet
{wang2012}, \citet{FFB12}], among others.

Before declaring that the folded concave penalty is superior to the
LASSO, we need to resolve a missing puzzle in the picture. Theoretical
properties of the folded concave penalization are established for a
theoretic local solution. However, we have to employ one of these local
minimization algorithms to find such a local optimal solution. It
remains to prove that the computed local solution is the desired
theoretic local solution to make the theory fully relevant. Many have
tried to address this issue [\citet{mcp}, \citet{fan2011},
\citet{zhangzhang2012}]. The basic idea there is to find
conditions under which the folded concave penalization actually has a
unique {sparse local} minimizer, and hence eliminate the problem of
multiple local solutions. Although this line of thoughts is natural and
logically intuitive, the imposed conditions for the unique {sparse
local} minimizer are too strong to be realistic.

In this paper, we offer a different approach to directly deal with the
multiple local solutions issue. We outline a general procedure based on
the LLA algorithm for solving a specific local solution of the folded
concave penalization, and then derive a lower bound on the probability
that this specific local solution exactly equals the oracle estimator.
This probability lower bound equals \mbox{$1-\delta_0-\delta_1-\delta_2$},
where $\delta_0$ corresponds to the exception probability of the
localizability of the underlying model, $\delta_1$ and $\delta_2$
represent the exception probability of the regularity of the oracle
estimator and they have nothing to do with any actual estimation
method. Explicit expressions of $\delta_0$, $\delta_1$ and $\delta
_2$ are given in Section~\ref{sec2}. Under weak regularity conditions, $\delta
_1$ and $\delta_2$ are very small. Thus, if $\delta_0$ goes to zero
then the computed solution is the oracle estimator with an overwhelming
probability. On the other hand, if $\delta_0$ cannot go to zero, then
it means that the model is extremely difficult to estimate no matter
how clever an estimator is. Thus, our theory suggests a
``bet-on-folded-concave-penalization'' principle: as long as there is a
reasonable initial estimator, our procedure can deliver an optimal
estimator using the folded concave penalization via the one-step LLA
implementation. Once the oracle estimator is obtained, the LLA
algorithm converges in the next iteration, namely, the oracle estimator
is a fixed point. Furthermore, we use four concrete examples to show
that exception probabilities $\delta_0$, $\delta_1$ and $\delta_2$
go to zero at a fast rate under the ultra-high-dimensional setting
where $\log(p)=O(n^{\eta})$ for $\eta\in(0,1)$.

Throughout this paper, the following notation is used. For $\mathbf{U}
=(u_{ij})_{k\times l}$, let $\|\mathbf{U}\|_{\min}=\min_{(i,j)}|u_{ij}|$ be
its minimum absolute value, and let $\lambda_{\min}(\mathbf{U})$ and
$\lambda_{\max}(\mathbf{U})$ be its smallest and largest
eigenvalues. We
introduce several matrix norms: the $\ell_1$ norm $\|\mathbf{U}\|
_{\ell
_1}=\max_j\sum_i|u_{ij}|$, the $\ell_2$ norm $\|\mathbf{U}\|_{\ell
_2}={\lambda^{1/2}_{\max}(\mathbf{U}'\mathbf{U})}$, the $\ell
_\infty$ norm $\|
\mathbf{U}\|_{\ell_\infty}=\max_i\sum_j|u_{ij}|$, the entrywise
$\ell
_{1}$ norm $\|\mathbf{U}\|_{1}=\sum_{(i,j)}|u_{ij}|$ and the
entrywise $\ell
_{\infty}$ norm $\|\mathbf{U}\|_{\max}=\max_{(i,j)}|u_{ij}|$.

\section{Main results}\label{sec2}

We begin with an abstract presentation of the sparse estimation
problem. Consider estimating a model based on $n$ \emph{i.i.d.}
observations. The target of estimation is ``parameter'' $\bolds{\beta
}^{\star
}=(\beta_1^{\star},\ldots,\beta_p^{\star})'$, that is, the model
is parameterized by $\bolds{\beta}^{\star}$. The dimension $p$ is larger
than the sample size $n$. In some problems, the target $\bolds{\beta
}^{\star
}$ can be a matrix (e.g., an inverse covariance matrix). In such cases,
it is understood that $(\beta_1^{\star},\ldots,\beta_p^{\star})'$
is the vectorization of the matrix~$\bolds{\beta}^{\star}$. Denote the
support set as $\mathcal{A}=\{j\dvtx \beta_j^{\star}\neq0\}$ and its
cardinality is $s=|\mathcal{A}|$. The sparsity assumption means that
$s \ll p$.

Suppose that our estimation scheme is to get a local minimizer of the
following folded concave penalized estimation problem:
%
\begin{equation}
\label{eq1} \min_{\bolds{\beta}}\ell_n(\bolds{
\beta})+P_{\lambda} \bigl(|\bolds{\beta}| \bigr)
\end{equation}
with $\ell_n(\bolds{\beta})$ is a convex loss and $P_{\lambda
}(|\bolds{\beta}
|)=\sum_jP_{\lambda}(|\beta_j|)$ is a folded concave penalty. In our
general theory, $\ell_n(\bolds{\beta})$ in (\ref{eq1}) does not
need to be
differentiable. The above formulation is a bit abstract but covers many
important statistical models. For example, $\ell_n(\bolds{\beta})$
can be
the squared error loss in penalized least squares, the check loss in
penalized quantile regression and the negative log-quasi-likelihood
function in penalized maximum quasi-likelihood.

An oracle knows the true support set, and defines the oracle estimator as
%
\begin{equation}
\label{eq2} \hat{\bolds{\beta}}{}^{\mathrm{oracle}} = \bigl(\hat{\bolds{
\beta}}{}^{\mathrm{oracle}}_{\mathcal{A}},\mathbf{0} \bigr) = \arg\min
_{\bolds{\beta}\dvtx \bolds{\beta}_{\mathcal{A}^c}=\mathbf{0}
}\ell_n(\bolds{\beta}).
\end{equation}
We assume that (\ref{eq2}) is regular such that the oracle solution is
unique, namely,
%
\begin{equation}
\label{eq3} \nabla_j \ell_n \bigl(\hat{\bolds{
\beta}}{}^{\mathrm{oracle}} \bigr) = 0\qquad\forall j \in\mathcal{A},
\end{equation}
where $\nabla_j$ denotes the subgradient with respect to the $j$th element of $\bolds{\beta}$. If the convex loss is
differentiable, the
subgradient is the usual gradient.
The oracle estimator is not a feasible estimator but it can be used as
a theoretic benchmark for other estimators to compare with. An
estimator is said to have the oracle property if it has the same
asymptotic distribution as the oracle estimator [\citet{scad},
\citet{fan2004}]. Moreover, an estimator is said to have the
strong oracle property if the estimator equals the oracle estimator
with overwhelming probability [\citet{fan2011}].

Throughout this paper, we also assume that the penalty $P_{\lambda
}(|t|)$ is a general folded concave penalty function defined on $t\in
(-\infty,\infty)$ satisfying:
\begin{longlist}[(iii)]
\item[(i)] $P_{\lambda}(t)$ is increasing and concave in $t\in
[0,\infty)$ with $P_{\lambda}(0)=0$;
\item[(ii)] $P_{\lambda}(t)$ is differentiable in $t\in(0,\infty)$
with $P'_{\lambda}(0):=P'_{\lambda}(0+)\ge a_1\lambda$;\vspace*{1pt}
\item[(iii)] $P'_{\lambda}(t)\ge a_1\lambda$ for $t\in(0,
a_2\lambda]$;\vspace*{1pt}
\item[(iv)] $P'_{\lambda}(t)=0$ for $t\in[a\lambda,\infty)$ with
the pre-specified constant $a>a_2$.
\end{longlist}
Where $a_1$ and $a_2$ are two fixed positive constants. The above
definition follows and extends previous works on the SCAD and the MCP
[\citet{scad}, \citet{mcp}, \citet{fan2011}]. The
derivative of the SCAD penalty is
\[
P'_{\lambda}(t)=\lambda I_{\{t\le\lambda\}}+
\frac{(a\lambda
-t)_+}{a-1}I_{\{t>\lambda\}}\qquad\mbox{for some }a>2
\]
and the derivative of the MCP is
$
P'_{\lambda}(t)=(\lambda-\frac{t}{a})_+$, for some $a>1$.
It is easy to see that $a_1=a_2=1$ for the SCAD, and $a_1=1-a^{-1}$,
$a_2=1$ for the MCP. The hard-thresholding penalty
$P_{\lambda}(t)=\lambda^2-(t-\lambda)^2I_{\{t<\lambda\}}$
[\citet{fan2001}]
is another special case of the general folded concave penalty with
$a=a_1=1$, $a_2=\frac{1}2$.

Numerical results in the literature show that the folded concave
penalty outperforms the $\ell_1$ penalty in terms of estimation
accuracy and selection consistency. To provide understanding of their
differences, it is important to show that the obtained solution of the
folded concave penalization has better theoretical properties than the
$\ell_1$-penalization. The technical difficulty here is to show that
the computed local solution is the local solution with proven
properties. \citet{mcp} and \citet{fan2011} proved the
restricted global optimality that the oracle estimator is the unique
global solution in the subspace $\mathbb{S}_s$, which is the union of
all $s$-dimensional coordinate subspaces in $\mathbb{R}^p$. Under
strong conditions, \citet{zhangzhang2012} proved that the global
solution leads to desirable recovery performance and corresponds to the
unique sparse local solution, and hence any algorithm finding a sparse
local solution will find the desired global solution. The fundamental
problem with these arguments is that in reality it is very rare that
the concave regularization actually has a unique sparse local solution,
which in turn implies that these strong conditions are too stringent to
hold in practice. Evidence is given in the simulation studies in
Section~\ref{sec4} where we show that the concave penalization has multiple
sparse local solutions.

We argue that, although the estimator is defined via the folded concave
penalization, we only care about properties of the computed estimator.
It is perfectly fine that the computed local solution is not the global
minimizer, as long as it has the desired properties. In this paper, we
directly analyze a specific estimator by the local linear approximation
(LLA) algorithm [\citet{onestep}]. The LLA algorithm takes
advantage of the special folded concave structure and utilizes the
majorization--minimization (MM) principle to turn a concave
regularization problem into a sequence of weighted $\ell_1$ penalized
problems. Within each LLA iteration, the local linear approximation is
the best convex majorization of the concave penalty function [see
Theorem~2 of \citet{onestep}]. Moreover, the MM principle has
provided theoretical guarantee to the convergence of the LLA algorithm
to a stationary point of the folded concave penalization. By analyzing
a logarithmic number of the LLA iterations, \citet{tzhang2010}
and \citet{zhanghuang2012} proved that the LLA algorithm helps
the folded concave penalization reduce the estimation error of the
LASSO in sparse linear and generalized linear regression.
In contrast, thresholding Lasso will not have the oracle property if
the irresponsentable condition does not hold: once some important
variables are missed in the Lasso fit, they can not be rescued by
thresholding. This is a very different operation from the LLA algorithm.

Here, we summarize the details of the LLA algorithm as in Algorithm~\ref{algo1}.

\begin{algorithm}[t]
\caption{The local linear approximation (LLA) algorithm}\label{algo1}
\begin{enumerate}
\item Initialize $\hat{\bolds{\beta}}{}^{(0)}=\hat{\bolds{\beta}}{}^{\mathrm{initial}}$
and compute the adaptive weight
\[
\hat{\mathbf{w}}{}^{(0)} = \bigl(\hat w{}^{(0)}_{1},
\ldots, \hat w^{(0)}_{p} \bigr)' =
\bigl(P'_{\lambda} \bigl(\bigl|\hat\beta^{(0)}_{1}\bigr|
\bigr),\ldots,P'_{\lambda} \bigl(\bigl|\hat\beta^{(0)}_{p}\bigr|
\bigr) \bigr)'.
\]
\item For $m=1,2,\ldots,$ repeat the LLA iteration till convergence
\begin{itemize}[(2.a)]
\item[(2.a)] Obtain $\hat{\bolds{\beta}}{}^{(m)}$ by solving the following
optimization problem
\[
\hat{\bolds{\beta}}{}^{(m)}=\min_{\bolds{\beta}}
\ell_n(\bolds{\beta})+\sum_j \hat
w^{(m-1)}_j \cdot|\beta_j|,
\]
\item[(2.b)] Update the adaptive weight vector $\hat{\mathbf{w}}{}^{(m)}$ with
$\hat w^{(m)}_{j} =P'_{\lambda}(|\hat\beta^{(m)}_{j}|)$.
\end{itemize}
\end{enumerate}
\end{algorithm}

In the following theorems, we provide the nonasymptotic analysis of the
LLA algorithm for obtaining the oracle estimator in the folded concave
penalized problem if it is initialized by some appropriate initial
estimator. In particular, the convex loss $\ell_n(\bolds{\beta})$ is not
required to be differentiable. To simplify notation, define $\nabla
\ell_n(\bolds{\beta})=(\nabla_1\ell_n(\bolds{\beta}),\ldots,\nabla_p\ell
_n(\bolds{\beta}))$ as the subgradient vector of $\ell_n(\bolds
{\beta})$. Denote
by $\mathcal{A}^c=\{j\dvtx \beta^\star_j=0\}$ the complement\vspace*{2pt} of the true
support set $\mathcal{A}$, and set $\nabla_{\mathcal{A}^c}\ell
_n(\bolds{\beta})=(\nabla_{j}\ell_n(\bolds{\beta})\dvtx j\in\mathcal
{A}^c)$ with
respect to $\mathcal{A}^c$.

\begin{theorem}\label{nonconvexopt1} Suppose the minimal signal
strength of $\bolds{\beta}^{\star}$ satisfies that
\begin{enumerate}[(A0)]
\item[(A0)] $\|\bolds{\beta}^{\star}_{\mathcal{A}}\|_{\min
}>(a+1)\lambda$.
\end{enumerate}
Consider the folded concave penalized problem with $P_{\lambda}(\cdot
)$ satisfying \textup{(i)--(iv)}. Let $a_0=\min\{1,a_2\}$. Under the event
\[
\mathcal{E}_1= \bigl\{\bigl\|\hat{\bolds{\beta}}{}^{\mathrm{initial}}-\bolds{
\beta}^{\star}\bigr\| _{\max}\le a_0 \lambda \bigr\} \cap
\bigl\{\bigl\Vert\nabla_{\mathcal{A}^c}\ell_n \bigl(\hat{\bolds{\beta}}{}^{\mathrm{oracle}} \bigr)\bigr\Vert_{\max} < a_1\lambda \bigr\},
\]
the LLA algorithm initialized by $\hat{\bolds{\beta}}{}^{\mathrm{initial}}$ finds
$\hat{\bolds{\beta}}{}^{\mathrm{oracle}}$ after one iteration.
\end{theorem}

Applying the union bound to $\mathcal{E}_1$, we easily get the
following corollary.

\begin{corollary}\label{co1}
With probability at least $1-\delta_0-\delta_1$, the LLA algorithm
initialized by $\hat{\bolds{\beta}}{}^{\mathrm{initial}}$ finds $\hat{\bolds{\beta}}{}^{\mathrm{oracle}}$ after one iteration, where
\[
\delta_0=\Pr \bigl(\bigl\|\hat{\bolds{\beta}}{}^{\mathrm{initial}}-\bolds{
\beta}^{\star}\bigr\| _{\max}>a_0\lambda \bigr)
\]
and
\[
\delta_1=\Pr \bigl(\bigl\Vert\nabla_{\mathcal{A}^c}\ell_n
\bigl(\hat{\bolds{\beta}}{}^{\mathrm{oracle}} \bigr)\bigr\Vert_{\max} \ge
a_1\lambda \bigr).
\]
\end{corollary}

\begin{rem}\label{rem1}
By its definition, $\delta_0$ represents the
localizability of the underlying model. To apply Theorem \ref{nonconvexopt1}, we need to have an appropriate initial estimator to
make $\delta_0$ go to zero as $n$ and $p$ diverge to infinity, namely
the underlying problem is localizable. In Section~\ref{sec3}, we will show by
concrete examples on how to find a good initial estimator to make the
problem localizable. $\delta_1$ represents the regularity behavior of
the oracle estimator, that is, its closeness to the true ``parameter''
measured by the score function. Note that
$\nabla_{\mathcal{A}^c}\ell_n(\bolds{\beta}^{\star})$ is concentrated
around zero. Thus, $\delta_1$ is usually small. In summary, Theorem
\ref{nonconvexopt1} and its corollary state that as long as the
problem is localizable and regular, we can find an oracle estimator by
using the one-step local linear approximation, which is a
generalization of the one-step estimation idea [\citet{onestep}]
to the high-dimensional setting.\looseness=1
\end{rem}

\begin{theorem}\label{nonconvexopt2}
Consider the folded concave penalized problem (\ref{eq1}) with
$P_{\lambda}(\cdot)$ satisfying \textup{(i)--(iv)}. Under the event
\[
\mathcal{E}_2= \bigl\{\bigl\|\nabla_{\mathcal{A}^c}\ell_n
\bigl(\hat{\bolds{\beta}}{}^{\mathrm{oracle}} \bigr)\bigr\|_{\max} <
a_1 \lambda \bigr\}\cap \bigl\{\bigl\| \hat{\bolds{\beta}}{}^{\mathrm{oracle}}_{\mathcal{A}}
\bigr\|_{\min} > a\lambda \bigr\},
\]
if $\hat{\bolds{\beta}}{}^{\mathrm{oracle}}$ is obtained, the LLA algorithm
will find
$\hat{\bolds{\beta}}{}^{\mathrm{oracle}}$ again in the next iteration, that
is, it
converges to $\hat{\bolds{\beta}}{}^{\mathrm{oracle}}$ in the next iteration
and is a
fixed point.
\end{theorem}

Now we combine Theorems \ref{nonconvexopt1} and \ref{nonconvexopt2}
to derive the nonasymptotic probability bound for the LLA algorithm to
exactly converge to the oracle estimator.

\begin{corollary}\label{nonconvexopt3}
Consider the folded concave penalized problem (\ref{eq1}) with
$P_{\lambda}(\cdot)$ satisfying \textup{(i)--(iv)}. Under assumption~\textup{(A0)}, the
LLA algorithm initialized by $\hat{\bolds{\beta}}{}^{\mathrm{initial}}$
converges to
$\hat{\bolds{\beta}}{}^{\mathrm{oracle}}$ after two iterations with
probability at
least $1-\delta_0-\delta_1-\delta_2$, where
\[
\delta_2=\Pr \bigl(\bigl\|\hat{\bolds{\beta}}{}^{\mathrm{oracle}}_{\mathcal
{A}}
\bigr\|_{\min
}\le a\lambda \bigr).
\]
\end{corollary}

\begin{rem}\label{rem2}
The localizable probability $1-\delta_0$ and
regularity probability $1-\delta_1$ have been defined before. $\delta
_2$ is a probability on the magnitude of the oracle estimator. Both
$\delta_1$ and $\delta_2$ are related to the regularity behavior of
the oracle estimator and will be referred to the oracle regularity condition.
Under assumption~\textup{(A0)}, it requires only the uniform convergence of
$\hat{\bolds{\beta}}{}^{\mathrm{oracle}}_{\mathcal{A}}$. Namely,
\[
\delta_2 \leq\Pr \bigl( \bigl\| \hat{\bolds{\beta}}{}^{\mathrm{oracle}}_{\mathcal{A}}
- \bolds{\beta}^\star_{\mathcal{A}}\bigr\|_{\max} > \lambda \bigr).
\]
Thus, we can regard $\delta_2$ as a direct measurement of the
closeness of the oracle estimator to the true ``parameter'' and is
usually small because of the small intrinsic dimensionality $s$. This
will indeed be shown in Section~\ref{sec3}.
\end{rem}

\begin{rem}\label{rem3}
The philosophy of our work follows the
well-known one-step estimation argument [\citet{bickel1975}] in
the maximum likelihood estimation. In some likelihood models, the
log-likelihood function is not concave. One of the local maximizers of
the log-likelihood is shown to be asymptotic efficient, but how to
compute that estimator is very challenging. \citet{bickel1975}
overcame this difficulty by focusing on a specially designed one-step
Newton--Raphson estimator initialized by a root-$n$ estimator. This
one-step estimator is asymptotically efficient, just like the
theoretical MLE. Note that Bickel's theory did not try to get the
global maximizer nor the theoretical local maximizer of the
log-likelihood, although the log-likelihood was used to construct the
explicit estimator. Our general theory follows this line of thinking.
Theorems \ref{nonconvexopt1}--\ref{nonconvexopt2} show how to
construct the explicit estimator that possesses the desired strong
oracle property. This is all we need to close the theoretical gap.
Following \citet{bickel1975}, we can just use the two-step LLA
solution and do not need to care about whether the LLA algorithm
converges. Of course, Theorem~\ref{nonconvexopt2} does offer a statistical convergence
proof of the LLA algorithm, which differs from its numeric convergence
argument. Moreover, Theorem~\ref{nonconvexopt2} also proves the statistical equivalence
between the two-step and fully converged LLA solutions. Thus, we
recommend using the two-step LLA solution as the folded concave
penalized estimator in applications.
\end{rem}

\section{Theoretical examples}\label{sec3}

This section outlines four classical examples to demonstrate
interesting and powerful applications of Theorems \ref{nonconvexopt1}
and \ref{nonconvexopt2}. We consider the linear regression, logistic
regression, precision matrix estimation and quantile regression.
We basically need to check the localizable condition and the regularity
condition for these problems.

\subsection{Sparse linear regression}\label{sec3.1}

The first example is the canonical problem of the folded concave
penalized least square estimation, that is,
%
\begin{equation}
\label{eq4} \min_{\bolds{\beta}}\frac{1}{2n}\|\mathbf{y}-
\mathbf{X}\bolds{\beta}\| _{\ell_2}^2+\sum
_jP_{\lambda} \bigl(|\beta_j| \bigr),
\end{equation}
where $\mathbf{y}\in\mathbb{R}^n$ and $\mathbf{X}=(\mathbf
{x}_1,\mathbf
{x}_2,\ldots,\mathbf{x}_n)'\in\mathbb{R}^{n\times p}$. Let
$\bolds{\beta}^{\star}$ be the true parameter vector in the linear
regression model $\mathbf{y}=\mathbf{X}\bolds{\beta}^{\star
}+\varepsilon
$, and the true
support set of $\bolds{\beta}^{\star}=(\beta^{\star}_j)_{1\le j\le
p}$ is
$\mathcal{A}=\{j\dvtx \beta^{\star}_j\neq0\}$. For the folded concave
penalized least square problem, the oracle estimator has an explicit
form of
\[
\hat{\bolds{\beta}}{}^{\mathrm{oracle}}= \bigl(\hat{\bolds{\beta }}{}^{\mathrm{oracle}}_{\mathcal
{A}},
\mathbf{0} \bigr)\qquad\mbox{with } \hat{\bolds{\beta}}{}^{\mathrm{oracle}}_{\mathcal{A}}=
\bigl(\mathbf{X}_{\mathcal
{A}}'\mathbf{X} _{\mathcal{A}}
\bigr)^{-1}\mathbf{X}_{\mathcal{A}}'\mathbf{y}
\]
and the Hessian matrix of $\ell_n(\bolds{\beta})$ is $n^{-1} \mathbf
{X}'\mathbf{X}$
regardless of $\bolds{\beta}$. Applying Theorems~\mbox{\ref{nonconvexopt1}--\ref{nonconvexopt2}}, we can derive the following theorem with explicit
upper bounds for $\delta_1$ and $\delta_2$, which depends only on the
behavior of the oracle estimator.

\begin{theorem}\label{lse}
Recall that $\delta_0=\Pr(\|\hat{\bolds{\beta}}{}^{\mathrm{initial}}-\bolds{\beta}
^{\star}\|_{\max}>a_0\lambda)$. Suppose
\begin{enumerate}[(A1)]
\item[(A1)] $\mathbf{y}=\mathbf{X}\bolds{\beta}^{\star
}+\varepsilon$
with $\varepsilon
=(\varepsilon_1,\ldots,\varepsilon_n)$ being i.i.d. sub-Gaussian
$(\sigma)$ for some fixed constant $\sigma>0$, that is,
$E[\exp(t\varepsilon_i^2)]\le\exp(\sigma^2t^2/2)$.
\end{enumerate}
The LLA algorithm initialized by $\hat{\bolds{\beta}}{}^{\mathrm{initial}}$ converges
to $\hat{\bolds{\beta}}{}^{\mathrm{oracle}}$ after two iterations with
probability at least
$1-\delta_0-\delta_1^{\mathrm{linear}}-\delta_2^{\mathrm{linear}}$, where
\[
\delta_1^{\mathrm{linear}}=2(p-s)\cdot\exp \biggl(-\frac{a_1^2n\lambda
^2}{2M\sigma^2}
\biggr)
\]
and
\[
\delta_2^{\mathrm{linear}}=2s\cdot\exp \biggl(-\frac{n \lambda_{\min
}}{2\sigma^2}\cdot
\bigl(\bigl\|\bolds{\beta}^{\star}_{\mathcal{A}}\bigr\|_{\min}-a\lambda
\bigr)^2 \biggr),
\]
where $\lambda_{\min} =\lambda_{\min}(\frac{1}n\mathbf
{X}'_{\mathcal
{A}}\mathbf{X}_{\mathcal{A}})$ and $M=\max_{j}\frac{1}n\|\mathbf
{x}_{(j)}\|_{\ell_2}^2$, which is usually 1 due to normalization, with
$\mathbf{x}_{(j)}=(x_{1j},\ldots,x_{nj})'$.
\end{theorem}

By Theorem~\ref{lse}, $\delta_1^{\mathrm{linear}}$ and $\delta_2^{\mathrm{linear}}$ go
to zero very quickly. Then it remains to bound $\delta_0$. To analyze
$\delta_0$, we should decide the initial estimator. Here, we use the
LASSO [\citet{lasso}] to initialize the LLA algorithm, which is
%
\begin{equation}
\label{eq5} \hat{\bolds{\beta}}{}^{\mathrm{lasso}}=\arg\min_{\bolds{\beta}}
\frac
{1}{2n}\|\mathbf{y}-\mathbf{X} \bolds{\beta}\|_{\ell_2}^2+
\lambda_{\mathrm{lasso}}\|\bolds{\beta}\|_{\ell_1}.
\end{equation}
Note that $\hat{\bolds{\beta}}{}^{\mathrm{lasso}}$ is the one-step LLA solution
initialized by zero. In order to bound $\hat{\bolds{\beta
}}{}^{\mathrm{lasso}}-\bolds{\beta}
^{\star}$, we invoke the following restricted eigenvalue condition:
\begin{enumerate}[(C1)]
\item[(C1)]
$
\kappa_{\mathrm{linear}}=\min_{\mathbf{u}\neq0\dvtx \|\mathbf{u}_{\mathcal
{A}^c}\|_{\ell_1}\le3\|\mathbf{u}_{\mathcal{A}}\|
_{\ell_1}}
\frac{\|\mathbf{X}\mathbf{u}\|^2_{\ell_2}}{n\|\mathbf{u}\|
^2_{\ell_2}}\in(0,\infty)$.
\end{enumerate}
This condition was studied in \citet{bickel2009}, \citet
{geer2009} and \citet{wainwright2012}. Under the assumptions of
\textup{(A1)} and~\textup{(C1)}, the LASSO yields the unique optimum $\hat{\bolds{\beta}}{}^{\mathrm{lasso}}$ satisfying
\[
\bigl\|\hat{\bolds{\beta}}{}^{\mathrm{lasso}}-\bolds{\beta}^{\star}
\bigr\|_{\ell
_2}\le\frac{4s^{1/2}\lambda_{\mathrm{lasso}}}{\kappa_{\mathrm{linear}}}
\]
with probability at least $1-2p\exp(-\frac{n\lambda
_{\mathrm{lasso}}^2}{2M\sigma^2})$. Thus, using this as an upper bound for $\|
\hat{\bolds{\beta}}{}^{\mathrm{lasso}}-\bolds{\beta}^{\star} \|_{\max}$,
it is easy for
us to obtain the following corollary.

\begin{corollary}\label{co3}
Under assumptions \textup{(A0)}, \textup{(A1)} and \textup{(C1)}, if we pick $\lambda\ge\frac
{4s^{1/2}\lambda_{\mathrm{lasso}}}{a_0\kappa_{\mathrm{linear}}}$, the LLA algorithm\vspace*{-2pt}
initialized by $\hat{\bolds{\beta}}{}^{\mathrm{lasso}}$ converges to
$\hat{\bolds{\beta}}{}^{\mathrm{oracle}}$ after two iterations with probability at least $1-2p\exp
(-\frac{n\lambda_{\mathrm{lasso}}^2}{2M\sigma^2})-\delta_1^{\mathrm{linear}}-\delta
_2^{\mathrm{linear}}$.
\end{corollary}

\begin{rem}\label{rem4}
Corollary \ref{co3} also suggests that sometimes it
is good to use zero to initialize the LLA algorithm.\vadjust{\goodbreak} If $\hat{\bolds
{\beta}}{}^{\mathrm{initial}}=\mathbf{0}$, the first LLA iteration gives a LASSO
estimator with $\lambda_{\mathrm{lasso}}=P'_{\lambda}(0)$. For both SCAD
and MCP, $P'_{\lambda}(0)=\lambda$.
If \mbox{$\lambda_{\mathrm{lasso}}=\lambda$} satisfies requirements in Corollary \ref{co3},
then after two more LLA iterations, the LLA algorithm converges to the
oracle estimator with high probability.
To be more specific, we have the following corollary.
\end{rem}

\begin{corollary}\label{co4}
Consider the SCAD or MCP penalized linear regression. Under assumptions
\textup{(A0)}, \textup{(A1)} and \textup{(C1)}, if $a_0\kappa_{\mathrm{linear}}\ge{4s^{1/2}}$, the LLA
algorithm \mbox{initialized} by zero converges to the oracle estimator after
three iterations with probability at least $1-2p\cdot\exp(-\frac
{n\lambda^2}{2M\sigma^2})-\delta_1^{\mathrm{linear}}-\delta_2^{\mathrm{linear}}$.
\end{corollary}

\begin{rem}\label{rem5}
The $s^{1/2}$ factor appears in Corollaries
\ref{co3}--\ref{co4} because $\|\hat{\bolds{\beta}}{}^{\mathrm{lasso}}-\bolds{\beta}^{\star
} \|_{\ell_2}$
is used to bound $\|\hat{\bolds{\beta}}{}^{\mathrm{lasso}}-\bolds{\beta
}^{\star} \|_{\max
}$. {It is possible to get rid of the $s^{1/2}$ factor by using the
$\ell_\infty$ loss of the LASSO in \citet{tzhang2009} and
\citet{ye-zhang-2010}. \citet{ye-zhang-2010} introduced the
cone invertability factor condition, that is,
\begin{enumerate}[(C1$'$)]
\item[(C1$'$)]
$
\zeta_{\mathrm{linear}}=\min_{\mathbf{u}\neq0\dvtx \|\mathbf{u}_{\mathcal
{A}^c}\|_{\ell_1}\le3\|\mathbf{u}_{\mathcal{A}}\|
_{\ell_1}}
\frac{\|\mathbf{X}'\mathbf{X}\mathbf{u}\|_{\max}}{n\|
\mathbf{u}\|
_{\max}}\in(0,\infty)$.
\end{enumerate}
Under assumptions \textup{(A1)} and \textup{(C1$'$)}, the LASSO yields
$\hat{\bolds{\beta}}{}^{\mathrm{lasso}}$ satisfying\break
$
\|\hat{\bolds{\beta}}{}^{\mathrm{lasso}}-\bolds{\beta}^{\star}\|_{{\max
}}\le
\frac{3\lambda_{\mathrm{lasso}}}{2\zeta_{\mathrm{linear}}}
$
with probability at least $1-2p\exp(-\frac{n\lambda
_{\mathrm{lasso}}^2}{8M\sigma^2})$.
}
\end{rem}

\begin{rem}\label{rem6}
Although we have considered using the LASSO
as the initial estimator, we can also use the Dantzig selector
[\citet{Dantzig07}] as the initial estimator, and the same
analysis can still go through under a similar restricted eigenvalue
condition as in \citet{bickel2009} {or a similar cone
invertability factor condition as in \citet{ye-zhang-2010}}.
\end{rem}

\subsection{Sparse logistic regression}\label{sec3.2}

The second example is the folded concave penalized logistic regression.
Assume that
\begin{enumerate}[(A2)]
\item[(A2)] the conditional distribution of $y_i$ given $\mathbf
{x}_i$ ($i=1,2,\ldots,n$) is a Bernoulli distribution with
$\Pr(y_i=1|\mathbf{x}_i,\bolds{\beta}^{\star})={\exp(\mathbf
{x}_i'\bolds{\beta}^{\star})}/(1+\exp(\mathbf{x}_i'\bolds{\beta
}^{\star}))$.
\end{enumerate}
Then the penalized logistic regression is given by
%
\begin{equation}
\label{eq6} \min_{\bolds{\beta}}\frac{1}{n}\sum
_i \bigl\{-y_i\mathbf{x}_i'
\bolds{\beta} +\psi \bigl(\mathbf{x}_i'\bolds{\beta}
\bigr) \bigr\}+\sum_jP_{\lambda
}\bigl(|\beta_j|\bigr)
\end{equation}
with the canonical link $\psi(t)=\log(1+\exp(t))$. This is a
canonical model for high-dimensional classification problems, and it is
a classical example of generalized linear models. The oracle estimator
is given by
\[
\hat{\bolds{\beta}}{}^{\mathrm{oracle}} = \bigl(\hat{\bolds{\beta
}}{}^{\mathrm{oracle}}_{\mathcal{A}}, \mathbf{0} \bigr) =\arg\min
_{\bolds{\beta}\dvtx \bolds{\beta}_{\mathcal
{A}^c}=\mathbf{0}
}\frac
{1}{n}\sum_i
\bigl\{-y_i\mathbf{x}_i' \bolds{\beta}+
\psi \bigl(\mathbf{x}_i'\bolds{\beta} \bigr) \bigr\}.
\]

For ease of presentation, we define
$
\bolds{\mu}(\bolds{\beta})=(\psi'(\mathbf{x}_1'\bolds{\beta
}),\ldots,\psi
'(\mathbf{x}_n'\bolds{\beta}))'
$
and
$
\bolds{\Sigma}(\bolds{\beta})=\operatorname{diag}\{\psi''(\mathbf
{x}_1'\bolds{\beta
}),\ldots,\psi''(\mathbf{x}_n'\bolds{\beta})\}$.
We also define three useful quantities:
$Q_1=\max_{j}\lambda_{\max}(\frac{1}n\mathbf{X}_{\mathcal
{A}}'\operatorname{diag}\{
|\mathbf{x}_{(j)}|\}\mathbf{X}_{\mathcal{A}})$, $Q_2=\|(\frac
{1}n\mathbf{X}
_{\mathcal{A}}'\bolds{\Sigma}(\bolds{\beta}^{\star})\mathbf
{X}_{\mathcal
{A}})^{-1}\|
_{\ell_\infty}$,
and
$Q_3=\|\mathbf{X}_{\mathcal{A}^c}'\bolds{\Sigma}(\bolds{\beta
}^{\star
})\mathbf{X}_{\mathcal
{A}}(\mathbf{X}_{\mathcal{A}}'\bolds{\Sigma}(\bolds{\beta}^{\star
})\mathbf{X}_{\mathcal
{A}})^{-1}\|_{\ell_\infty}$,
where $\operatorname{diag}\{|\mathbf{x}_{(j)}|\}$ is a diagonal
matrix with
elements $\{|x_{ij}|\}_{i=1}^n$.

%
\begin{theorem}\label{logistic}
Recall that $\delta_0=\Pr(\|\hat{\bolds{\beta}}{}^{\mathrm{initial}}-
\bolds{\beta}^{\star}\|_{\max}>a_0\lambda)$. Under assumption~\textup{(A2)},
the LLA algorithm initialized by $\hat{\bolds{\beta}}{}^{\mathrm{initial}}$ converges
to $\hat{\bolds{\beta}}{}^{\mathrm{oracle}}$ after two iterations with
probability at least
$1-\delta_0-\delta_1^{\mathrm{logit}}-\delta_2^{\mathrm{logit}}$, where
\begin{eqnarray*}
\delta_1^{\mathrm{logit}} &=& 2s\cdot\exp \biggl(-\frac{n}{M}
\min \biggl\{\frac
{2}{Q_1^2Q_2^4s^2},\frac{a_1^2\lambda^2}{2(1+2Q_3)^2} \biggr\} \biggr)
\\
&&{} +2(p-s)\cdot \exp \biggl(-\frac{a_1^2n\lambda^2}{2M} \biggr)
\end{eqnarray*}
with $M = \max_{j} n^{-1} \| \mathbf{x}_{(j)}\|_{\ell_2}^2$ and
\[
\delta_2^{\mathrm{logit}}= 2s\cdot\exp \biggl(-\frac{n}{MQ_2^2}\min
\biggl\{\frac
{2}{Q_1^2Q_2^2s^2},\frac{1}2 \bigl(\bigl\|\bolds{
\beta}_{\mathcal
{A}}{}^{\star
} \bigr\| _{\min}-a\lambda
\bigr)^2 \biggr\} \biggr).
\]
\end{theorem}

Under fairly weak assumptions, $\delta_1^{\mathrm{logit}}$ and $\delta
_2^{\mathrm{logit}}$ go to zero very quickly. The remaining challenge is to
bound $\delta_0$. We consider using the $\ell_1$-penalized maximum
likelihood estimator as the initial estimator, that is,
\[
\hat{\bolds{\beta}}{}^{\mathrm{lasso}} =\arg\min_{\bolds{\beta}}
\frac{1}{n}\sum_i \bigl\{-y_i
\mathbf{x}_i'\bolds{\beta}+\psi \bigl(
\mathbf{x}_i'\bolds{\beta} \bigr) \bigr\} +\lambda
_{\mathrm{lasso}}\|\bolds{\beta}\|_{\ell_1}.
\]

\begin{theorem}\label{logisticlasso}
Let $m=\max_{(i,j)}|x_{ij}|$. Under assumption \textup{(A2)} and
\begin{enumerate}[(C2)]
\item[(C2)] $ \kappa_{\mathrm{logit}}=\min_{\mathbf{u}\neq
\mathbf{0}\dvtx \|\mathbf{u}_{\mathcal{A}^c}\|_{\ell_1}\le3\|
\mathbf{u}_{\mathcal{A}}\|_{\ell_1}}\frac{\mathbf{u}'\nabla^2\ell
^{\mathrm{logit}}_n(\bolds{\beta}^{\star})\mathbf{u}}{\mathbf{u}'\mathbf{u}}\in
(0,\infty)$,
\end{enumerate}
if $\lambda_{\mathrm{lasso}}\le\frac{\kappa_{\mathrm{logit}}}{20ms}$, with
probability at least
$
1-2p\cdot\exp(-\frac{n}{2M}\lambda_{\mathrm{lasso}}^2)$,
we have
\[
\bigl\|\hat{\bolds{\beta}}{}^{\mathrm{lasso}}-\bolds{\beta}^{\star}
\bigr\|_{\ell
_2}\le5\kappa_{\mathrm{logit}}^{-1}s^{1/2}
\lambda_{\mathrm{lasso}}.
\]
\end{theorem}

In light of Theorem \ref{logisticlasso}, we can obtain the following
corollary.

\begin{corollary}\label{co5}
Under assumptions \textup{(A0)}, \textup{(A2)} and \textup{(C2)}, if we pick $\lambda\ge\frac
{5s^{1/2}\lambda_{\mathrm{lasso}}}{a_0\kappa_{\mathrm{logit}}}$, the LLA algorithm\vspace*{-2pt}
initialized by $\hat{\bolds{\beta}}{}^{\mathrm{lasso}}$ converges to
$\hat{\bolds{\beta}}{}^{\mathrm{oracle}}$ after two iterations with probability at least $1-2p\exp
(-\frac{n}{2M}\lambda_{\mathrm{lasso}}^2)-\delta_1^{\mathrm{logit}}-\delta_2^{\mathrm{logit}}$.
\end{corollary}

Again we can use zero to initialize the LLA algorithm and do three LLA
iterations, because the first LLA iteration gives a $\ell_1$ penalized
logistic regression with $\lambda_{\mathrm{lasso}}=P'_{\lambda}(0)$ which
equals $\lambda$ for both SCAD and MCP.

\begin{corollary}\label{co6}
Consider the SCAD/MCP penalized logistic regression. Under assumptions
\textup{(A0)}, \textup{(A2)} and \textup{(C2)}, if $a_0\kappa_{\mathrm{logit}}\ge5s^{1/2}$ holds, the LLA
algorithm initialized by zero converges to the oracle estimator after
three iterations with probability at least $1-2p\exp(-\frac
{n}{2M}\lambda^2)-\delta_1^{\mathrm{logit}}-\delta_2^{\mathrm{logit}}$.
\end{corollary}

\begin{rem}\label{rem7}
The $s^{1/2}$ factor appears in Corollaries
\ref{co5}--\ref{co6} because $\|\hat{\bolds{\beta}}{}^{\mathrm{lasso}}-\bolds{\beta}^{\star
} \|_{\ell_2}$
is used to bound $\|\hat{\bolds{\beta}}{}^{\mathrm{lasso}}-\bolds{\beta
}^{\star} \|_{\max
}$. To remove the $s^{1/2}$ factor, we can use the general
invertability factor condition [\citet{zhanghuang2012}] to obtain
the $\ell_\infty$ loss of the LASSO. For space consideration, details
are omitted.
\end{rem}

\subsection{Sparse precision matrix estimation}\label{sec3.3}

The third example is the folded concave penalized Gaussian
quasi-likelihood estimator for the sparse precision matrix estimation
problem, that is,
%
\begin{equation}
\label{eq7} \min_{\bolds{\Theta}\succ0}-\log\det(\bolds{\Theta})+\langle
\bolds{\Theta},\widehat{\bolds{\Sigma}}_n\rangle+\sum
_{(j,k)\dvtx j\neq k}P_{\lambda}\bigl(|\theta_{jk}|\bigr)
\end{equation}
with the sample covariance matrix $\widehat{\bolds{\Sigma}}_n=(\hat
\sigma
^n_{ij})_{q\times q}$. Under the Gaussian assumption, the sparse
precision matrix is translated into a sparse Gaussian graphical model.
In this example, the target ``parameter'' is the true precision matrix
$\bolds{\Theta}^{\star}=(\theta^{\star}_{jk})_{q\times q}$ with the
support set $\mathcal{A}=\{(j,k)\dvtx \theta_{jk}^{\star}\neq0\}$. Due
to the symmetric structure of $\bolds{\Theta}$, the dimension of the target
``parameter'' is $p=q(q+1)/2$, and the cardinality of $\mathcal{A}$ is
$s=\#\{(j,k)\dvtx j\le k, \theta_{jk}^{\star}\neq0\}$. Moreover, we
denote the maximal degree of $\bolds{\Theta}^{\star}$ as $d=\max_j\#
\{
k\dvtx \theta_{jk}^{\star}\neq0\}$.

In the sparse precision matrix estimation, the oracle estimator is
given by
\[
\widehat{\bolds{\Theta}}{}^{\mathrm{oracle}} =\arg\min_{\bolds{\Theta}\succ\mathbf
{0}\dvtx \bolds{\Theta
}_{\mathcal
{A}^c}=\mathbf{0}
}-\log
\det(\bolds{\Theta})+\langle\bolds{\Theta},\widehat{\bolds{\Sigma}}_n
\rangle.
\]
The Hessian matrix of $\ell_n(\bolds{\Theta})$ is $\mathbf
{H}^{\star}=\bolds{\Sigma}
^{\star}\otimes\bolds{\Sigma}^{\star}$. To simplify notation,
we let
\[
K_1=\bigl\|\bolds{\Sigma}^{\star}\bigr\|_{\ell_\infty}, \qquad
K_2=\bigl\| \bigl(\mathbf{H}^{\star
}_{\mathcal{A}\mathcal{A}}
\bigr)^{-1}\bigr\|_{\ell_\infty}\quad\mbox{and} \quad K_3=\bigl\|
\mathbf{H}^{\star}_{\mathcal{A}^c\mathcal{A}} \bigl(\mathbf{H}^{\star
}_{\mathcal{A}\mathcal{A}}
\bigr)^{-1}\bigr\|_{\ell_\infty}.
\]

In the next theorem, we derive the explicit bounds for $\delta_1$ and
$\delta_2$ under the Gaussian assumption. Similar results hold under
the exponential/polynomial tail condition in \citet{clime}. Under
the normality, we cite a large deviation result [\citet
{saulis1991}, \citet{bickel2008a}]:
%
\begin{equation}
\label{eq8} \Pr \bigl(|\hat\sigma^n_{ij}-
\sigma^\star_{ij}|\ge\nu \bigr)\le C_0\exp
\bigl(-c_0n\nu^2 \bigr)
\end{equation}
for any $\nu$ such that $|\nu|\le\nu_0$, where $\nu_0$, $c_0$ and
$C_0$ depend on $\max_{i}\sigma^\star_{ii}$ only.

\begin{theorem}\label{glasso}
Let $\delta_0^{G}=\Pr(\|\widehat{\bolds{\Theta}}{}^{\mathrm{initial}}-\bolds
{\Theta}
^{\star}\|_{\max}>a_0\lambda)$. Assume that
\begin{enumerate}[(A0$'$)]
\item[(A$0'$)] $\|\bolds{\Theta}^{\star}_{\mathcal{A}}\|
_{\min
}>(a+1)\lambda$, and
\item[(A3)] $\mathbf{x}_1,\ldots,\mathbf{x}_n$ are
i.i.d. Gaussian samples with the true covariance $\bolds
{\Sigma}
^{\star}$.
\end{enumerate}
The LLA algorithm initialized by $\widehat{\bolds{\Theta}}{}^{\mathrm{initial}}$
converges to $\widehat{\bolds{\Theta}}{}^{\mathrm{oracle}}$ after two iterations with
probability at least $1-\delta_0^G-\delta_1^G-\delta_2^G$, where
\begin{eqnarray*}
\delta_1^G &=& C_0s\cdot\exp \biggl(-
\frac{c_0}{4}n\cdot\min \biggl\{\frac
{a_1^2\lambda^2}{(2K_3+1)^2},\frac{1}{9K_1^{2}K_2^{2}d^2},
\frac
{1}{9K_1^{6}K_2^{4}d^2} \biggr\} \biggr)
\\
&&{}+ C_0(p-s)\cdot\exp \biggl(-\frac{c_0a_1^2}{4}n
\lambda^2 \biggr)
\end{eqnarray*}
and
\[
\delta_2^G = C_0s\cdot\exp \biggl(-
\frac{c_0n}{4K_2^2}\cdot{\min \biggl\{\frac{1}{9K_1^{2}d^{2}},\frac
{1}{9K_1^{6}K_2^{2}d^2}, \bigl(
\bigl\|\bolds{\Theta}^{\star}_{\mathcal{A}}\bigr\|_{\min}-a\lambda
\bigr)^2 \biggr\} } \biggr).
\]
\end{theorem}

Theorem~\ref{glasso} shows that both $\delta_1^G$ and $\delta_2^G$
go to zero very quickly. Now we only need to deal with $\delta_0^G$.
To initialize the LLA algorithm, we consider using the constrained
$\ell_1$ minimization estimator (CLIME) by \citet{clime}, that is,
\[
\widehat{\bolds{\Theta}}{}^{\mathrm{clime}}=\arg\min_{\bolds{\Theta}}\|
\bolds{\Theta}\|_1\qquad\mbox{subject to }\|\widehat{\bolds{
\Sigma}}_n\bolds{\Theta}-\mathbf{I}\| _{\max}\le
\lambda_{\mathrm{clime}}.
\]

Define $L=\|\bolds{\Theta}^{\star}\|_{\ell_1}$. As discussed in
\citet
{clime}, it is reasonable to assume that $L$ is upper bounded by a
constant or $L$ is some slowly diverging quantity, because $\bolds
{\Theta}
^{\star}$ has a few nonzero entries in each row. We combine the
concentration bound (\ref{eq8}) and the same line of proof as in
\citet{clime} to show that with probability at least $1-C_0p\cdot
\exp(-\frac{c_0n}{L^2}\lambda^2_{\mathrm{clime}})$, we have
\[
\bigl\|\widehat{\bolds{\Theta}}{}^{\mathrm{clime}}-\bolds{\Theta}^{\star}
\bigr\|_{\max
}\le4L\lambda_{\mathrm{clime}}.
\]

Thus we have the following corollary.
%
\begin{corollary}\label{co7}
Under assumptions \textup{(A0$'$)} and \textup{(A3)}, if $\lambda\ge\frac
{4L}{a_0}\lambda_{\mathrm{clime}}$, the LLA algorithm initialized by $\widehat
{\bolds{\Theta}}{}^{\mathrm{clime}}$ converges to $\widehat{\bolds{\Theta
}}{}^{\mathrm{oracle}}$ after two
iterations with probability at least $1-C_0p\exp(-\frac
{c_0n}{L^2}\lambda^2_{\mathrm{clime}})-\delta_1^{G}-\delta_2^{G}$.
\end{corollary}

In the current literature, the $\ell_1$ penalized likelihood estimator
GLASSO [\citet{yuan2007}] is perhaps the most popular estimator
for sparse precision matrix estimation. However, it requires a strong
irrepresentable condition [\citet{ravikumar2008}] stating that
there exists a fixed constant $\gamma_G\in(0,1)$ such that
$
\|\mathbf{H}^{\star}_{\mathcal{A}^c\mathcal{A}}(\mathbf{H}^{\star
}_{\mathcal
{A}\mathcal{A}})^{-1}\|_{\ell_\infty}\le\gamma_G$.
This condition is very restrictive. If we replace the $\ell_1$ penalty
with a folded concave penalty, it is interesting to see that we can
obtain the oracle precision matrix estimator by using CLIME as the
initial estimator in the LLA algorithm without requiring any strong
structure assumption such as the irrepresentable condition.

\subsection{Sparse quantile regression}\label{sec3.4}

The fourth example is the folded concave penalized quantile regression.
Quantile regression [\citet{koenker2005}] has wide applications
in statistics and econometrics. Recently, the spare quantile regression
has received much attention [\citet{cqr}, \citet{zhu2008},
\citet{wu2009}, \citet{belloni2011},
\citet{wang2012}, \citet{FFB12}]. We consider estimating
the conditional $\tau$ quantile under
\begin{enumerate}[(A4)]
\item[(A4)] $\mathbf{y}=\mathbf{X}\bolds{\beta}^{\star
}+\varepsilon$
with $\varepsilon
=(\varepsilon_1,\ldots,\varepsilon_n)$ being the independent errors
satisfying $\Pr(\varepsilon_i\le0) = \tau$ for some fixed constant
$\tau\in(0,1)$. Let $f_i(\cdot)$ be the density function of
$\varepsilon_i$, and define $F_i(\cdot)$ as its distribution function.
\end{enumerate}

Denote by $\rho_{\tau}(u)=u\cdot(\tau-I_{\{u\le0\}})$ the check
loss function [\citet{koenker1978}]. The folded concave penalized
quantile regression is given by
%
\begin{equation}
\label{eq9} \min_{\bolds{\beta}}\frac{1}{n}\sum
_i \rho_{\tau} \bigl(y_i-\mathbf
{x}_i'\bolds{\beta} \bigr) + \sum
_jP_{\lambda}\bigl(|\beta_j|\bigr).
\end{equation}
The oracle estimator of the sparse quantile regression is given by
\[
\hat{\bolds{\beta}}{}^{\mathrm{oracle}} = \bigl(\hat{\bolds{\beta
}}{}^{\mathrm{oracle}}_{\mathcal{A}}, \mathbf{0} \bigr) =\arg\min
_{\bolds{\beta}\dvtx \bolds{\beta}_{\mathcal
{A}^c}=\mathbf{0}
}\frac
{1}{n}\sum_i
\rho_{\tau} \bigl(y_i- \mathbf{x}_i'
\bolds{\beta} \bigr).
\]

Note that the check loss $\rho_{\tau}(\cdot)$ is convex but
nondifferentiable. Thus, we need to handle the subgradient $\nabla\ell
_n(\bolds{\beta})=(\nabla_1\ell_n(\bolds{\beta}),\ldots, \nabla
_p\ell_n(\bolds{\beta}
))$, where
\[
\nabla_j\ell_n(\bolds{\beta})= \frac{1}n\sum
_i x_{ij}\cdot \bigl((1-\tau)
I_{\{y_i - \mathbf
{x}_i' \bolds{\beta}{\le} 0\}}-z_jI_{\{y_i - \mathbf{x}_i' \bolds
{\beta}=0\}
}-\tau I_{\{y_i - \mathbf{x}_i' \bolds{\beta}>0\}} \bigr)
\]
with $z_j\in[\tau-1, \tau]$ is the subgradient of $\rho_{\tau}(u)$
when $u=0$. To simplify notation, we let $M_{\mathcal{A}}=\max_{i}\frac
{1}s\|\mathbf{x}_{i\mathcal{A}}\|_{\ell_2}^2$, and
$m_{\mathcal{A}^c}=\max_{(i,j)\dvtx j\in\mathcal{A}^c}|x_{ij}|$.

\begin{theorem}\label{quantile}
Recall that $\delta_0=\Pr(\|\hat{\bolds{\beta}}{}^{\mathrm{initial}}-\bolds{\beta}
^{\star}\|_{\max}>a_0\lambda)$. Suppose
\begin{enumerate}[(C3)]
\item[(C3)] there exist constants $u_0>0$ and $0<f_{\min}\le f_{\max}
<\infty$ such that for any $u$ satisfying $|u|\le u_0$, $f_{\min} \le
\min_if_i(u)\le\max_if_i(u)\le f_{\max}$.
\end{enumerate}
%
If $\lambda \gg 1/n$ such that $\log p = o(n\lambda^2)$, $(M_{\mathcal
{A}}s)^{1/2}(\|\bolds{\beta}^\star_{\mathcal{A}}\|_{\min}-a\lambda
)\le
u_0$, and $m_{\mathcal{A}^c}M_{\mathcal{A}}s=o(\frac{n^{1/2}\lambda
}{\log^{1/2} n})$, the LLA algorithm\vspace*{-2pt} initialized by $\hat{\bolds{\beta}}{}^{\mathrm{initial}}$ converges to $\hat{\bolds{\beta}}{}^{\mathrm{oracle}}$ after two
iterations with probability at least $1-\delta_0-\delta_1^{Q}-\delta
_2^{Q}$, where
$
\delta_1^{Q}
=
4n^{-1/2}+C_1 (p-s)\cdot\exp(-\frac{a_1n\lambda
}{104m_{\mathcal{A}^c}})+2(p-s)\cdot\exp(-\frac{a_1^2n\lambda
^2}{32m_{\mathcal{A}^c}^2})
$,\vspace*{-2pt} and
$
\delta_2^{Q}
=
4\exp(-\frac{\lambda_{\min}^2 f_{\min}^2}{72M_{\mathcal
{A}}}\cdot\frac{n}s(\|\bolds{\beta}^\star_{\mathcal{A}}\|_{\min
}-a\lambda)^2 )
$
with $\lambda_{\min} =\lambda_{\min}(\frac{1}n\mathbf
{X}'_{\mathcal
{A}}\mathbf{X}_{\mathcal{A}})$ and $C_1>0$ that does not depend on
$n$, $p$
or $s$.
\end{theorem}

Under fairly weak assumptions, both $\delta_1^{Q}$ and $\delta_2^{Q}$
go to zero very quickly. Next, we only need to bound $\delta_0$. We
consider using the $\ell_1$-penalized quantile regression as the
initial estimator in the LLA algorithm, that is,
\[
\hat{\bolds{\beta}}{}^{\mathrm{lasso}} =\arg\min_{\bolds{\beta}}
\frac{1}{n}\sum_i \rho_{\tau
}
\bigl(y_i-\mathbf{x}_i'\bolds{\beta}
\bigr)+\lambda_{\mathrm{lasso}}\|\bolds{\beta}\| _{\ell_1}.
\]
To bound $\delta_0$, we use the estimation bound for $\hat{\bolds{\beta}}{}^{\mathrm{lasso}}$ by \citet{belloni2011} and summarize the main result in
the following lemma.

\begin{lemma} \label{belloni2011}
Under assumption \textup{(A4)} and the assumption of Theorem 2 in \citet
{belloni2011}, which implies $\gamma\rightarrow0$, for any $A>1$ and
$p^{-1}\le\alpha\rightarrow0$, if $\lambda_{\mathrm{lasso}}$ satisfies (3.8)
of \citet{belloni2011}, with probability at least $1-\alpha
-4\gamma-3p^{-A^2}$, $\hat{\bolds{\beta}}{}^{\mathrm{lasso}}$ satisfies
\[
\bigl\|\hat{\bolds{\beta}}{}^{\mathrm{lasso}}-\bolds{\beta}^{\star}
\bigr\|_{\ell
_2}\le C_2\sqrt{\frac{s\log p}{n}},
\]
where $C_2>0$ is a fixed constant that does not depend on $s$, $p$ and $n$.
\end{lemma}

Thus we have the following corollary.

\begin{corollary}\label{co8}
Under the assumptions of Lemma \ref{belloni2011}, for any $A>1$ and
$p^{-1}\le\alpha\rightarrow0$, if $\lambda$ such that $\lambda\ge
\frac{C_2}{a_0}\sqrt{{s\log p}/{n}}$ and $\lambda$ also\vspace*{-1pt} satisfies
the conditions in Theorem \ref{quantile}, the LLA algorithm
initialized by $\hat{\bolds{\beta}}{}^{\mathrm{lasso}}$ converges to
$\hat{\bolds{\beta}}{}^{\mathrm{oracle}}$ after two iterations with probability at least $1-\alpha
-4\gamma-3p^{-A^2}-\break \delta_1^{Q}-\delta_2^{Q}$.
\end{corollary}

\section{Simulation studies}\label{sec4}

In this section, we use simulation to examine the finite sample
properties of the folded concave penalization for solving four
classical problems. We fixed $a=3.7$ in the SCAD and $a=2$ in the MCP
as suggested in \citet{scad} and \citet{mcp}, respectively.

\subsection{Sparse regression models}\label{sec4.1} We first considered sparse
linear, logistic and quantile regression models. In all examples,
we simulated $n$ training data and $n$ validation data and generated
$\mathbf{x}\sim N_p(0,\bolds{\Sigma})$ with $\bolds{\Sigma}
=(0.5^{|i-j|})_{p\times p}$.

\begin{mo}[(Sparse linear regression)]\label{mo1} Set $n=100$ and $p=1000$.
The response $y=\mathbf{x}'\bolds{\beta}^\star+\varepsilon$, where
$\bolds{\beta}^\star=(3, 1.5, 0, 0, 2, 0_{p-5})$ and $\varepsilon
\sim N(0,1)$.
\end{mo}
The validation error of a generic estimator $\hat{\bolds{\beta}}$ for
Model~\ref{mo1} is defined as
\[
\sum_{i \in\mathrm{validation}} \bigl(y_i-
\mathbf{x}'_i \hat{\bolds{\beta}} \bigr)^2.
\]

\begin{mo}[(Sparse logistic regression)]\label{mo2}
Set $n=200$ and $p=1000$.
The response $y$ follows a Bernoulli distribution with success
probability as\break  $\frac{\exp(\mathbf{x}'\bolds{\beta}^\star
)}{1+\exp
(\mathbf{x}'\bolds{\beta}^\star)}$,
where $\bolds{\beta}^\star$ is constructed by randomly choosing $10$
elements in $\bolds{\beta}^\star$ as $t_1s_1,\ldots,t_{10}s_{10}$ and
setting the other $p-10$ elements as zero, where $t_{j}$'s are
independently drawn from $\operatorname{Unif}(1,2)$, and $s_{j}$'s are
independent Bernoulli samples with $\Pr(s_{j}=1)=\Pr(s_{j}=-1)=0.5$.
\end{mo}
The validation error of a generic estimator $\hat{\bolds{\beta}}$ for
Model~\ref{mo2} is defined as
\[
\sum_{i \in\mathrm{validation}} \bigl(-y_i
\mathbf{x}'_i \hat{\bolds{\beta}}+\log \bigl(1+\exp
\bigl(\mathbf{x}'_i \hat{\bolds{\beta}} \bigr) \bigr)
\bigr).
\]

\begin{mo}[(Sparse quantile regression)]\label{mo3}
Set $n=100$ and $p=400$.
The response $y=\mathbf{x}'\bolds{\beta}^\star+\varepsilon$ where
$\bolds{\beta}^\star$
is constructed in the same way as in Model~\ref{mo2}, and $\varepsilon$
follows the standard Cauchy distribution.
\end{mo}
We considered $\tau=0.3,0.5$ in the simulation. The validation error
of a generic estimator $\hat{\bolds{\beta}}$ for Model~\ref{mo3} is
defined as
\[
\sum_{i \in\mathrm{validation}} \rho_{\tau}
\bigl(y_i- \mathbf{x}_i'\hat{\bolds{
\beta}} \bigr).
\]
We include $\ell_1$ penalization in the simulation study, and computed
the $\ell_1$ penalized linear/logistic regression and quantile
regression by R packages \emph{glmnet} and \emph{quantreg}. We
implemented three local solutions of SCAD/MCP. The first\vadjust{\goodbreak} one, denoted
by SCAD-cd/MCP-cd, was the fully convergent solutions computed by
coordinate descent. The second one, denoted by SCAD-lla0/MCP-lla0, was
computed by the LLA algorithm initialized by zero. The third one,
denoted by SCAD-lla$^\star$/MCP-lla$^\star$, was computed by the LLA
algorithm initialized by the tuned LASSO estimator. SCAD-lla$^\star
$/MCP-lla$^\star$ is the fully iterative LLA solution designed
according to the theoretical analysis in Sections~\ref{sec3.1}, \ref{sec3.2} and \ref{sec3.4}. We
also computed the three-step LLA solution of SCAD-lla0/MCP-lla0,
denoted by SCAD-3slla0/MCP-3slla0, and the two-step LLA solution of
SCAD-lla$^\star$/MCP-lla$^\star$, denoted by SCAD-2slla$^\star
$/MCP-2slla$^\star$. Especially, SCAD-2slla$^\star$/MCP-2slla$^\star
$ is the recommended two-step LLA solution.
For each competitor, its penalization parameter was chosen by
minimizing the validation error.

We conducted 100 independent runs for each model. Estimation accuracy
is measured by the average $\ell_1$ loss $\|\hat{\bolds{\beta}}-\bolds
{\beta}
^\star\|_{\ell_1}$ and $\ell_2$ loss $\|\hat{\bolds{\beta}}-\bolds{\beta}
^\star\|_{\ell_2}$, and selection accuracy is evaluated by the
average counts of false positive and false negative. The simulation
results are summarized in Table~\ref{simulationm1m2m3}.

\begin{table}
\tabcolsep=5pt
\caption{Numerical comparison of LASSO, SCAD $\&$ MCP in Models \protect\ref{mo1}--\protect\ref{mo3}.
Estimation accuracy is measured by the $\ell_1$ loss and the $\ell_2$
loss, and selection accuracy is measured by counts of false negative (\#FN) or false positive (\#FP). Each metric is averaged over $100$
replications with its standard error shown in the parenthesis}\label{simulationm1m2m3}
{\fontsize{9pt}{9.85pt}\selectfont{\begin{tabular*}{\tablewidth}{@{\extracolsep{\fill}}@{}lcccc@{\qquad}cccc@{}}
\hline
\textbf{Method}
& \multicolumn{1}{c}{$\bolds{\ell_1}$ \textbf{loss}} & \multicolumn{1}{c}{$\bolds{\ell_2}$ \textbf{loss}} & \multicolumn{1}{c}{\textbf{\#FP}} & \multicolumn{1}{c}{\textbf{\#FN}}
& \multicolumn{1}{c}{$\bolds{\ell_1}$ \textbf{loss}} & \multicolumn{1}{c}{$\bolds{\ell_2}$ \textbf{loss}} & \multicolumn{1}{c}{\textbf{\#FP}} & \multicolumn{1}{c@{}}{\textbf{\#FN}}
\\
\hline
& \multicolumn{4}{c}{Model \ref{mo1} (linear regression)} & \multicolumn{4}{c}{Model \ref{mo2} (logistic regression)}\\
{LASSO}
&1.20 &0.46 &14.68 &0 &15.08 &3.62&55.92 &0.59 \\
&(0.05)&(0.01)&(0.74)&(0)&(0.06)&(0.02) &(0.93)&(0.04)\\[2pt]
{SCAD-cd}
&0.34 &0.20 &2.22 &0 &9.12 &2.40&27.72 &0.58 \\
&(0.02)&(0.01) &(0.40)&(0)&(0.15)&(0.04)&(0.46)&(0.04)\\[2pt]
{SCAD-3slla0}
&0.29 &0.20 &0 &0 &6.79 &2.44&0.90 &2.72 \\
&(0.01)&(0.01) &(0) &(0) &(0.13)&(0.03) &(0.06)&(0.08)\\[2pt]
{SCAD-lla0}
&0.29 &0.20 &0 &0 &6.42 &2.34&0.79 &2.73 \\
&(0.01)&(0.01) &(0) &(0) &(0.13)&(0.03) &(0.06)&(0.08)\\[2pt]
{SCAD-2slla$^\star$}
&0.30 &0.20 &0 &0 &6.65 &2.41&0.76 &2.58 \\
&(0.01)&(0.01) &(0) &(0) &(0.15)&(0.04)&(0.08)&(0.08)\\[2pt]
{SCAD-lla$^{\star}$}
&0.29 &0.19 &0 &0 &6.41 &2.33&0.74 &2.74 \\
&(0.01)&(0.01) &(0) &(0) &(0.14)&(0.04)&(0.06)&(0.09)\\[2pt]
{MCP-cd}
&0.31 &0.20 &0.75 &0 &6.97 &2.30&3.62 &1.46 \\
&(0.02)&(0.01) &(0.14)&(0)&(0.16)&(0.05)&(0.14)&(0.08)\\[2pt]
{MPC-3slla0}
&0.30 &0.20 &0 &0 &7.10 &2.52&0.94 &2.86\\
&(0.02)&(0.01)&(0) &(0) &(0.14)&(0.04)&(0.09)&(0.09)\\[2pt]
{MPC-lla0}
&0.29 &0.20 &0 &0 &6.88 &2.45&1.11 &2.81 \\
&(0.02)&(0.01)&(0) &(0) &(0.14)&(0.04)&(0.09)&(0.09)\\[2pt]
{MCP-2slla$^\star$}
&0.29 &0.19 &0 &0 &6.79 &2.43&0.96 &2.49 \\
&(0.02)&(0.01)&(0) &(0)&(0.14)&(0.04)&(0.08)&(0.09)\\[2pt]
{MCP-lla$^{\star}$}
&0.28 &0.19 &0 &0 &6.30 &2.30&0.78 &2.64 \\
&(0.02)&(0.01)&(0) &(0) &(0.14)&(0.04)&(0.07)&(0.08)\\[5pt]
& \multicolumn{8}{c@{}}{Model \ref{mo3} (quantile regression)}\\
& \multicolumn{4}{c}{$\tau=0.3$} & \multicolumn{4}{c@{}}{$\tau=0.5$}\\
{LASSO}
&14.33 &2.92 &39.31 &1.03 &13.09 &2.62 &41.42 &0.61 \\
&(0.35) &(0.07) &(1.29) &(0.14) &(0.33) &(0.06) &(1.18) &(0.09) \\[2pt]
{SCAD-3slla0}
&9.08 &2.31 &22.79 &1.27 &6.58 &1.65 &22.43 &0.62\\
&(0.43) &(0.10) &(1.02) &(0.15) &(0.38) &(0.09) &(1.03) &(0.11) \\[2pt]
{SCAD-lla0}
&7.70 &2.20 &16.08 &1.66 &4.47 &1.37 &13.48 &0.68 \\
&(0.46) &(0.12) &(0.94) &(0.20) &(0.31) &(0.09) &(0.72) &(0.12) \\[2pt]
{SCAD-2slla$^\star$}
&7.43 &2.13 &13.26 &1.43 &4.80 &1.50 &11.43 &0.74 \\
&(0.45) &(0.10) &(1.02) &(0.15) &(0.29) &(0.08) &(0.75) &(0.11) \\[2pt]
{SCAD-lla$^\star$}
&5.92 &1.93 &8.89 &1.63 &3.96 &1.27 &10.18 &0.69 \\
&(0.37) &(0.11) &(0.68) &(0.19) &(0.39) &(0.08) &(0.74) &(0.11) \\[2pt]
{MCP-3slla0}
&10.09 &2.71 &19.34 &1.58 &7.44 &1.93 &17.54 &0.88 \\
&(0.44) &(0.09) &(1.10) &(0.17) &(0.44) &(0.10) &(0.77) &(0.14)\\[2pt]
{MCP-lla0}
&9.86 &2.69 &11.63 &2.18 &5.79 &1.70 &9.45 &1.04 \\
&(0.53) &(0.12) &(0.95) &(0.21) &(0.44) &(0.10) &(0.77) &(0.14) \\[2pt]
{MCP-2slla$^\star$}
&6.13 &2.15 &5.10 &1.75 &4.48 &1.54 &3.53 &1.03 \\
&(0.44) &(0.10) &(0.54) &(0.16) &(0.29) &(0.09) &(0.47) &(0.13) \\[2pt]
{MCP-lla$^\star$}
&5.95 &2.05 &2.92 &1.91 &3.88 &1.39 &2.04 &1.00 \\
&(0.42) &(0.11) &(0.50) &(0.18) &(0.29) &(0.09) &(0.29) &(0.13)\\
\hline
\end{tabular*}}}
\end{table}

We can draw the following main conclusions:
\begin{longlist}[(3)]
\item[(1)] The local solutions solved by coordinate descent and LLA are
different. LLA using different initial values are technically different
algorithms, and the corresponding fully\vadjust{\goodbreak} converged solutions are
different, also. This message clearly shows that the concave
regularization problem does have multiple local minimizers and multiple
sparse local minimizers.

\item[(2)] SCAD-2slla$^\star$/MCP-2slla$^\star$ are recommended based
on our philosophy and theory (see Remark~\ref{rem3}). SCAD-2slla$^\star
$/MCP-2slla$^\star$ is asymptotically equivalent to SCAD-lla$^\star
$/MCP-lla$^\star$, but two-step solutions are cheaper to compute than
fully converged ones. We do expect to see they differ with finite
sample size, but the difference is ignorable as in Table~\ref{simulationm1m2m3}.

\item[(3)] We included SCAD-3slla0/MCP-3slla0 because of Corollaries \ref{co4}
and \ref{co6} that justify the use of zero as a good initial value in the LLA
algorithm under some extra conditions. The simulation results show that
zero can be a good initial value but it is not the best choice one
would try.
\end{longlist}

\subsection{Sparse Gaussian graphical model}\label{sec4.2}

We drew $n=100$ training data and $n$ validation data from
$N_{q=100}(\mathbf{0},\bolds{\Sigma}^\star)$ with a sparse precision matrix~$\bolds{\Theta}^\star$.

\begin{mo}\label{mo4}
$\bolds{\Theta}^\star$ is a tridiagonal matrix by setting
$\bolds{\Sigma}^\star=(\sigma^\star_{ij})_{q\times q}$ as an AR(1)
covariance matrix with $\sigma^\star_{ij}=\exp(-|s_i-s_j|)$ for
$s_1<\cdots<s_q$, which draws $s_q-s_{q-1}, s_{q-1}-s_{q-2},\ldots,s_2-s_1$ independently from $\operatorname{Unif}(0.5,1)$.
\end{mo}

\begin{mo}\label{mo5}
$\bolds{\Theta}^\star=\mathbf{U}'_{q\times
q}\mathbf{U}_{q\times q}+\mathbf{I}
_{q\times q}$ where $\mathbf{U}=(u_{ij})_{q\times q}$ has zero
diagonals and
$100$ nonzero off-diagonal entries. The nonzero entries are generated
by $u_{ij}=t_{ij}s_{ij}$ where $t_{ij}$'s are independently drawn from
$\operatorname{Unif}(1,2)$, and $s_{ij}$'s are independent Bernoulli variables
with $\Pr(s_{ij}=\pm1)=0.5$.
\end{mo}
The validation error of a generic estimator $\widehat{\bolds{\Theta}}$ for
Models \ref{mo4}--\ref{mo5} is defined as
\[
-\log\det(\widehat{\bolds{\Theta}})+ \bigl\langle\widehat{\bolds{\Theta}},
\widehat{\bolds{\Sigma}}_n^{\mathrm{validation}} \bigr\rangle.
\]

We computed the GLASSO and the CLIME by the R packages \emph{glasso}
and \emph{clime}. We computed two local solutions of the SCAD/MCP
penalized \mbox{estimator} denoted by GSCAD/GMCP. The first one, denoted
by\break
\mbox{GSCAD-lla0/GMCP-lla0}, used $\operatorname{diag}(\widehat\Sigma^{-1}_{jj})$ as
the initial solution in the LLA algorithm. The second one, denoted by
GSCAD-lla$^\star$/GMCP-lla$^\star$, used the tuned CLIME to
initialize the LLA algorithm. GSCAD-lla$^\star$/GMCP-lla$^\star$ was
designed according to the theoretical analysis in Section~\ref{sec3.3}. We
computed the three-step LLA solution of \mbox{GSCAD-lla0/GMCP-lla0}, denoted
by~GSCAD-3slla0/GMCP-3slla0, and the two-step LLA solution of\break
\mbox{GSCAD-lla$^\star$/GMCP-lla$^\star$}, denoted by GSCAD-2slla$^\star
$/GMCP-2slla$^\star$.
For each competitor, its penalization parameter was chosen by
minimizing the validation error.

For each model, we conducted 100 independent runs. Estimation accuracy
is measured by the average Operator norm loss $\|\widehat{\bolds{\Theta}}
-\bolds{\Theta}\|_{\ell_2}$ and Frobenius norm loss $\|\widehat{\bolds
{\Theta}}
-\bolds{\Theta}\|_F$, and selection accuracy is evaluated by the average
counts of false positive and false negative. The simulation results are
summarized in Table~\ref{simulationm4m5}. The main conclusions are
the same as those in spare regression models. First, the fully
converged LLA solutions are different with different initial values, so
it is impractical to try to prove that this problem has a unique
minimizer. Second, with the CLIME as the initial value, the two-step
LLA solutions perform
as well as the fully converged LLA solutions.

\begin{table}
\tabcolsep=5.9pt
\caption{Numerical comparison of GLASSO, CLIME, GSCAD $and$ GMCP in
Model~\protect\ref{mo5}. Estimation accuracy is measured by the Operator norm ($\|\cdot
\|_{\ell_2}$) and the Frobenius norm ($\|\cdot\|_F$), and selection
accuracy is measured by counts of false negative (\#FN) or false
positive (\#FP). Each metric is averaged over $100$ replications with
its standard error in the parenthesis}\label{simulationm4m5}
\begin{tabular*}{\tablewidth}{@{\extracolsep{\fill}}lccccc@{\qquad}cccc@{}}
\hline
\textbf{Method} & $\bolds{\|\cdot\|_{\ell_2}}$ & $\bolds{\|\cdot\|_F}$ & \textbf{\#FP} & \textbf{\#FN}&
$\bolds{\|\cdot\|_{\ell_2}}$ & $\bolds{\|\cdot\|_F}$ & \textbf{\#FP} & \textbf{\#FN}\\
\hline
& \multicolumn{4}{c}{Model \ref{mo4}}& \multicolumn{4}{c@{}}{Model \ref{mo5}}\\
{GLASSO}
&1.45 &6.12 &743.56 &1.34 &11.63 &25.45 &236.76&56.16\\
&(0.01)&(0.02)&(10.75)&(0.17) &(0.02)&(0.03)&(5.19)&(0.52)\\[2pt]
{CLIME}
&1.40 &5.89 &741.16 &2.42 &8.56 &18.40 &323.04&12.26\\
&(0.01)&(0.03)&(12.80)&(0.24) &(0.05)&(0.08)&(7.22)&(0.38)\\[2pt]
{GSCAD-3slla0}
&1.20 &4.59 &659.04&1.78 &10.84 &22.05 &225.36&54.98\\
&(0.02)&(0.03)&(9.41)&(0.20) &(0.05)&(0.12)&(4.92)&(0.58)\\[2pt]
{GSCAD-lla0}
&1.16 &4.42 &641.82&1.96 &10.73 &20.68 &228.70&54.54\\
&(0.02)&(0.03)&(9.41)&(0.20) &(0.05)&(0.12)&(4.92)&(0.58)\\[2pt]
{GSCAD-2slla$^\star$}
&1.20 &4.60 &660.84&1.74 &6.49 &13.69 &203.52&28.78\\
&(0.02)&(0.03)&(9.39)&(0.19) &(0.13)&(0.15)&(5.27)&(0.57)\\[2pt]
{GSCAD-lla$^\star$}
&1.16 &4.42 &635.49&1.94 &6.42 &13.36 &196.60&30.02\\
&(0.02)&(0.03)&(9.39)&(0.19) &(0.13)&(0.15)&(5.27)&(0.57)\\[2pt]
{GMCP-3slla0}
&1.57 &4.62 &349.36&3.56 &10.41 &20.40 &201.06&54.12\\
&(0.04)&(0.04)&(7.03)&(0.21) &(0.07)&(0.17)&(4.26)&(0.68)\\[2pt]
{GMCP-lla0}
&1.53 &4.56 &291.04&6.45 &10.34 &19.20 &200.37&52.24\\
&(0.04)&(0.04)&(5.12)&(0.32) &(0.07)&(0.12)&(4.26)&(0.60)\\[2pt]
{GMCP-2slla$^\star$}
&1.40 &4.37 &290.10&3.96 &5.98 &12.74 &68.14 &23.84\\
&(0.04)&(0.04)&(4.71)&(0.20) &(0.17)&(0.17)&(3.42)&(0.61)\\[2pt]
{GMCP-lla$^\star$}
&1.39 &4.31 &229.87&6.29 &5.80 &12.72 &44.79 &25.18\\
&(0.03)&(0.04)&(4.56)&(0.33) &(0.15)&(0.17)&(3.12)&(0.53)\\
\hline
\end{tabular*}
\end{table}

\section{Technical proofs}\label{sec5}

\subsection{Proof of Theorem \texorpdfstring{\protect\ref{nonconvexopt1}}{1}}\label{sec5.1}
Define
$
\hat{\bolds{\beta}}{}^{(0)}= \hat{\bolds{\beta}}{}^{\mathrm{initial}}
$. Under the event
$
\{\| \hat{\bolds{\beta}}{}^{(0)}-\bolds{\beta}^{\star}\|_{\max}
\le a_0\lambda\}
$, due to assumption \textup{(A0)}, we have
$
|\hat\beta^{(0)}_{j}| \le\| \hat{\bolds{\beta}}{}^{(0)}-\bolds
{\beta}^{\star}\|
_{\max}\le a_0\lambda\le a_2\lambda
$
for $j \in{\mathcal{A}^c}$, and
$
|\hat\beta^{(0)}_j| \ge\|\bolds{\beta}^\star_{\mathcal{A}}\|
_{\min}-\|
\hat{\bolds{\beta}}{}^{(0)}-\bolds{\beta}^{\star}\|_{\max
}>a\lambda
$
for $j \in{\mathcal{A}}$. By property (iv), $P_\lambda'(|\hat\beta
^{(0)}_j|) = 0$ for all $j \in\mathcal{A}$. Thus $\hat{\bolds{\beta}}{}^{(1)}$ is the solution to the problem
%
\begin{equation}
\label{eq10} \hat{\bolds{\beta}}{}^{(1)} =\arg\min_{\bolds{\beta}}
\ell_n(\bolds{\beta})+\sum_{j \in
{\mathcal{A}}^c}
P'_{\lambda} \bigl(\bigl|\hat\beta^{(0)}_j\bigr|
\bigr)\cdot|\beta_{j}|.
\end{equation}
By properties (ii) and (iii), $P'_{\lambda}(|\hat\beta^{(0)}_j|)\ge
a_1\lambda$ holds for $j\in{\mathcal{A}}^c$. We now show that
$\hat{\bolds{\beta}}{}^{\mathrm{oracle}}$ is the unique global solution to
(\ref
{eq10}) under the additional condition $\{\|\nabla_{\mathcal
{A}^c}\ell_n(\hat{\bolds{\beta}}{}^{\mathrm{oracle}})\|_{\max} < a_1\lambda
\}$. To
see this, note that by convexity, we have
%
\begin{eqnarray}
\label{eq11} \ell_n(\bolds{\beta}) & \geq& \ell_n
\bigl(\hat{\bolds{\beta}}{}^{\mathrm{oracle}} \bigr) + \sum
_j \nabla_j \ell_n \bigl(\hat{
\bolds{ \beta}}{}^{\mathrm{oracle}} \bigr) \bigl(\beta_j - \hat
\beta^{\mathrm{oracle}}_j \bigr)
\nonumber\\[-8pt]\\[-8pt]
& = & \ell_n \bigl(\hat{\bolds{\beta}}{}^{\mathrm{oracle}} \bigr) + \sum
_{j \in
\mathcal{A}^c} \nabla_j \ell_n \bigl(
\hat{\bolds{\beta}}{}^{\mathrm{oracle}} \bigr) \bigl(\beta_j - \hat\beta
^{\mathrm{oracle}}_j \bigr),\nonumber
\end{eqnarray}
where the last equality used (\ref{eq3}). By (\ref{eq11}) and
$\hat{\bolds{\beta}}{}^{\mathrm{oracle}}_{\mathcal{A}^c}=\mathbf{0}$, for any
$\bolds{\beta}$
we have
\begin{eqnarray*}
&& \biggl\{ \ell_n(\bolds{\beta})+\sum_{j \in{\mathcal{A}}^c}
P'_{\lambda
} \bigl(\bigl|\hat\beta^{(0)}_j\bigr|
\bigr)|\beta_{j}| \biggr\} - \biggl\{\ell_n \bigl(\hat{
\bolds{\beta}}{}^{\mathrm{oracle}} \bigr)+\sum_{j \in
{\mathcal
{A}}^c}
P'_{\lambda} \bigl(\bigl|\hat\beta^{(0)}_j\bigr|
\bigr)\bigl|\hat\beta^{\mathrm{oracle}}_{j}\bigr| \biggr\}
\\
&&\qquad \ge \sum_{j \in\mathcal{A}^c} \bigl\{P'_{\lambda}
\bigl(\bigl|\hat\beta^{(0)}_j\bigr| \bigr)-\nabla_j
\ell_n \bigl(\hat{\bolds{\beta}}{}^{\mathrm{oracle}} \bigr)\cdot\operatorname{sign}(\beta_j) \bigr\} \cdot|\beta_{j}|
\\
&&\qquad \ge \sum_{j \in\mathcal{A}^c} \bigl\{ a_1\lambda-
\nabla_j\ell_n \bigl(\hat{\bolds{\beta}}{}^{\mathrm{oracle}}
\bigr)\cdot\operatorname{sign}(\beta_j) \bigr\} \cdot|\beta_{j}|
\\
&&\qquad \ge0.
\end{eqnarray*}
The strict inequality holds unless $\beta_j = 0$, $\forall j \in
\mathcal{A}^c$. This together with the uniqueness of the solution to
(\ref{eq2}) concludes that $\hat{\bolds{\beta}}{}^{\mathrm{oracle}}$ is the unique
solution to (\ref{eq10}). Hence, $\hat{\bolds{\beta}}{}^{(1)}=\hat{\bolds
{\beta}}{}^{\mathrm{oracle}}$, which completes the proof of Theorem \ref{nonconvexopt1}.

\subsection{Proof of Theorem \texorpdfstring{\protect\ref{nonconvexopt2}}{2}}\label{sec5.2}
Given that the LLA algorithm finds $\hat{\bolds{\beta}}{}^{\mathrm{oracle}}$
at the
current iteration, we denote $\hat{\bolds{\beta}}$ as the solution
to the
convex optimization problem in the next iteration of the LLA algorithm.
Using $\hat{\bolds{\beta}}{}^{\mathrm{oracle}}_{\mathcal{A}^c}=\mathbf{0}$ and
$P'_{\lambda}(|\hat\beta^{\mathrm{oracle}}_j|) = 0$ for $j \in\mathcal{A}$
under the event $\{\|\hat{\bolds{\beta}}{}^{\mathrm{oracle}}_{\mathcal{A}} \|
_{\min
} > a\lambda\}$, we have
%
\begin{equation}
\label{eq12} \hat{\bolds{\beta}}=\arg\min_{\bolds{\beta}}
\ell_n(\bolds{\beta})+\sum_{j \in
{\mathcal{A}}^c} \gamma
\cdot|\beta_{j}|,
\end{equation}
where $\gamma= P'_{\lambda}(0) \geq a_1\lambda$. This problem is
very similar to (\ref{eq10}). We can follow the proof of Theorem \ref
{nonconvexopt1} to show that under the additional condition $\{\|
\nabla_{\mathcal{A}^c}\ell_n(\hat{\bolds{\beta}}{}^{\mathrm{oracle}})\|
_{\max} <
a_1\lambda\}$, $\hat{\bolds{\beta}}{}^{\mathrm{oracle}}$ is the unique
solution to
(\ref{eq12}). Hence, the LLA algorithm converges, which completes the
proof of Theorem \ref{nonconvexopt2}.

\subsection{Proof of Theorem \texorpdfstring{\protect\ref{lse}}{3}}\label{sec5.3}
Let $\mathbf{H}_{\mathcal{A}}=\mathbf{X}_{\mathcal{A}}(\mathbf
{X}_{\mathcal{A}}'\mathbf{X}
_{\mathcal{A}})^{-1}\mathbf{X}_{\mathcal{A}}'$. Since $\mathbf
{y}=\break \mathbf
{X}_{\mathcal
{A}}\bolds{\beta}^{\star}_{\mathcal{A}}+\varepsilon$, we have
$
\nabla_{\mathcal{A}^c}\ell_n(\hat{\bolds{\beta}}{}^{\mathrm{oracle}})=
\frac{1}n\mathbf{X}_{\mathcal{A}^c}'(\mathbf{I}_{n\times n}-\mathbf
{H}_{\mathcal
{A}})\varepsilon.
$
By the Chernoff bound, we have
\begin{eqnarray*}
\delta_1 
&\le&\sum
_{j\in\mathcal{A}^c}\Pr \bigl(\bigl|\mathbf{x}_{(j)}'(
\mathbf{I} _{n\times n}-\mathbf{H}_{\mathcal{A}})\varepsilon\bigr|>
a_1n\lambda \bigr)
\\
&\le&2\sum_{j\in\mathcal{A}^c}\exp \biggl(-\frac{a_1^2n^2\lambda
^2}{2\sigma^2\cdot\|\mathbf{x}_{(j)}'(\mathbf{I}_{n\times
n}-\mathbf{H}
_{\mathcal{A}})\|_{\ell_2}^2}
\biggr).
\end{eqnarray*}
Since
$
\|\mathbf{x}_{(j)}'(\mathbf{I}_{n\times n}-\mathbf{H}_{\mathcal
{A}})\|_{\ell
_2}^2=\mathbf{x}_{(j)}'(\mathbf{I}_{n\times n}-\mathbf{H}_{\mathcal
{A}})\mathbf{x}_{(j)}\le nM$,
we conclude that
\[
\delta_1 \leq2(p-s)\exp \biggl(-\frac{a_1^2n\lambda^2}{2M\sigma^2} \biggr).
\]

Now\vspace*{1pt} we bound $\delta_2$. Note that
$
\hat{\bolds{\beta}}{}^{\mathrm{oracle}}_{\mathcal{A}}
=\bolds{\beta}^{\star}_{\mathcal{A}}+(\mathbf{X}_{\mathcal
{A}}'\mathbf{X}_{\mathcal
{A}})^{-1}\mathbf{X}_{\mathcal{A}}'\varepsilon$,
and then
$
\|\hat{\bolds{\beta}}{}^{\mathrm{oracle}}_{\mathcal{A}} \|_{\min}
\ge\|\bolds{\beta}^{\star}_{\mathcal{A}}\|_{\min}-\|(\mathbf
{X}_{\mathcal
{A}}'\mathbf{X}_{\mathcal{A}})^{-1}\mathbf{X}_{\mathcal
{A}}'\varepsilon\|_{\max}
$. Thus, we have
%
\begin{equation}
\label{eq13} \delta_2 \le\Pr \bigl(\bigl\| \bigl(\mathbf{X}_{\mathcal{A}}'
\mathbf{X}_{\mathcal
{A}} \bigr)^{-1}\mathbf{X}_{\mathcal
{A}}'
\varepsilon\bigr\|_{\max}\ge\bigl\|\bolds{\beta}^{\star}_{\mathcal
{A}}\bigr\|
_{\min}-a\lambda \bigr).
\end{equation}
It remains to derive an explicit bound for (\ref{eq13}). To simplify
notation, we let
$
(\mathbf{X}_{\mathcal{A}}'\mathbf{X}_{\mathcal{A}})^{-1}\mathbf
{X}_{\mathcal
{A}}'=(\mathbf{u}_1,\mathbf{u}_2,\ldots,\mathbf{u}_s)'$,
with $\mathbf{u}_j = \mathbf{X}_{\mathcal{A}} (\mathbf
{X}_{\mathcal{A}}'\mathbf{X}
_{\mathcal{A}})^{-1} \mathbf e_j$, where $\mathbf e_j$ is the
unit vector with $j$th element 1. It is obvious to see that
$
\|\mathbf{u}_j\|_{\ell_2}^2= \mathbf e_j' (\mathbf
{X}_{\mathcal
{A}}'\mathbf{X}_{\mathcal{A}})^{-1} \mathbf e_j' \leq( n \lambda
_{\min})^{-1}$.
By the Chernoff bound, we have
\begin{eqnarray*}
\delta_2 &\le&\Pr \bigl(\bigl\| \bigl(\mathbf{X}_{\mathcal{A}}'
\mathbf{X}_{\mathcal
{A}} \bigr)^{-1}\mathbf{X} _{\mathcal{A}}'
\varepsilon\bigr\|_{\max}\ge\bigl\|\bolds{\beta}^{\star
}_{\mathcal
{A}}
\bigr\|_{\min}-a\lambda \bigr)
\\
&\le& 2\sum_{j=1}^s\exp \biggl(-
\frac{(\|\bolds{\beta}^{\star
}_{\mathcal
{A}}\|_{\min}-a\lambda)^2}{2\sigma^2\|\mathbf{u}_{j}\|_{\ell
_2}^2} \biggr)
\\
&\le& 2s\exp \biggl(-\frac{n\lambda_{\min}}{2\sigma^2} \bigl(\bigl\|\bolds {\beta}
^{\star}_{\mathcal{A}}\bigr\|_{\min}-a\lambda \bigr)^2
\biggr).
\end{eqnarray*}
Thus, we complete the proof of Theorem \ref{lse}.

\subsection{Proof of Theorem \texorpdfstring{\protect\ref{logistic}}{4}}\label{sec5.4}

A translation of (\ref{eq3}) into our setting becomes
%
\begin{equation}
\label{eq14} \mathbf{X}_{\mathcal{A}}'\bolds{\mu} \bigl(\hat{
\bolds{\beta}}{}^{\mathrm{oracle}} \bigr)=\mathbf{X}_{\mathcal
{A}}'
\mathbf{y}.
\end{equation}
We now use this to bound $\delta_2$.

Define a map $F\dvtx \mathbb{B}(r)\subset\mathbb{R}^{p}\rightarrow
\mathbb{R}^{p}$ satisfying
$
F(\bolds{\Delta})= ((F_{\mathcal{A}}(\bolds{\Delta}_{\mathcal
{A}}))',\mathbf{0}
' )'
$
with
$
F_{\mathcal{A}}(\bolds{\Delta}_{\mathcal{A}})=(\mathbf
{X}_{\mathcal{A}}'\bolds{\Sigma}
(\bolds{\beta}^{\star})\mathbf{X}_{\mathcal{A}})^{-1}\cdot\mathbf
{X}_{\mathcal
{A}}'(\mathbf{y}-\bolds{\mu}(\bolds{\beta}^{\star}+\bolds{\Delta
}))+\bolds{\Delta}_{\mathcal{A}}
$
and the convex compact set
$
\mathbb{B}(r)=\{\bolds{\Delta}\in\mathbb{R}^p\dvtx \|\bolds{\Delta
}_{\mathcal{A}}\|
_{\max}\le r,\bolds{\Delta}_{\mathcal{A}^c}=\mathbf{0}\}
$
with $r=2Q_2\cdot\break \|\frac{1}n\mathbf{X}_{\mathcal{A}}'(\bolds{\mu
}(\bolds{\beta}^{\star
})-\mathbf{y})\|_{\max}$. Our aim is to show
%
\begin{equation}
\label{eq15} F \bigl(\mathbb{B}(r) \bigr)\subset\mathbb{B}(r),
\end{equation}
when
%
\begin{equation}
\label{eq16} \biggl\|\frac{1}n\mathbf{X}_{\mathcal{A}}'
\bigl(\bolds{\mu} \bigl(\bolds{\beta}^{\star} \bigr)-\mathbf{y} \bigr)
\biggr\|_{\max
}\le\frac{1}{Q_1Q_2^2s}.
\end{equation}

If\vspace*{1pt} (\ref{eq15}) holds, by the Brouwer's fixed-point theorem, there
always exists a fixed point $\widehat{\bolds{\Delta}}\in\mathbb
{B}(r)$ such
that $F(\widehat{\bolds{\Delta}})=\widehat{\bolds{\Delta}}$. It
immediately follows that
$
\mathbf{X}_{\mathcal{A}}'\mathbf{y}=\mathbf{X}_{\mathcal
{A}}'\bolds{\mu
}(\bolds{\beta}^{\star
}+\widehat{\bolds{\Delta}})
$
and
$
\widehat{\bolds{\Delta}}_{\mathcal{A}^c}=\mathbf{0}$,
which implies $\bolds{\beta}^{\star}+\widehat{\bolds{\Delta}}=\hat
{\bolds{\beta}}{}^{\mathrm{oracle}}$ by uniqueness of the solution to (\ref{eq14}). Thus,
%
\begin{equation}
\label{eq17} \bigl\|\hat{\bolds{\beta}}{}^{\mathrm{oracle}}-\bolds{\beta}^{\star}
\bigr\|_{\max
}=\|\widehat{\bolds{\Delta}}\|_{\max}\le r.
\end{equation}
If further
\[
\biggl\|\frac{1}n\mathbf{X}_{\mathcal{A}}' \bigl(\bolds{\mu}
\bigl(\bolds{\beta}^{\star} \bigr)-\mathbf{y} \bigr)\biggr\|_{\max
}\le
\frac{1}{2Q_2} \bigl(\bigl\|\bolds{\beta}_{\mathcal{A}}{}^{\star}
\bigr\|_{\min
}-a\lambda \bigr),
\]
we have
$
r \leq\|\bolds{\beta}_{\mathcal{A}}{}^{\star}\|_{\min}-a\lambda$,
and then
$
\|
\hat{\bolds{\beta}}{}^{\mathrm{oracle}}_{\mathcal{A}}
\|_{\min}
\ge
a\lambda$.
Therefore, we have
\[
\delta_2 \le\Pr \biggl(\biggl\|\frac{1}n\mathbf{X}_{\mathcal{A}}'
\bigl(\bolds{\mu} \bigl(\bolds{\beta}^{\star} \bigr)-\mathbf{y} \bigr)
\biggr\|_{\max}>\min \biggl\{\frac{1}{Q_1Q_2^2s},\frac{1}{2Q_2} \bigl(\bigl\|
\bolds{\beta}_{\mathcal{A}}{}^{\star}\bigr\|_{\min}-a\lambda \bigr)
\biggr\} \biggr).
\]
By the Hoeffding's bound in Proposition 4(a) of \citet{fan2011},
we have
\[
\delta_2 \leq2s\cdot\exp \biggl(-\frac{n}{MQ_2^2}\cdot\min \biggl
\{\frac
{2}{Q_1^2Q_2^2s^2},\frac{1}2 \bigl(\bigl\|\bolds{\beta}_{\mathcal{A}}{}^{\star
}
\bigr\| _{\min}-a\lambda \bigr)^2 \biggr\} \biggr).
\]

We now derive (\ref{eq15}). Using its Taylor expansion around $\bolds
{\Delta}
=\mathbf{0}$, we have
\[
\mathbf{X}_{\mathcal{A}}'\bolds{\mu} \bigl(\bolds{
\beta}^{\star
}+\bolds{\Delta} \bigr) = \mathbf{X}_{\mathcal{A}}'
\bolds{\mu} \bigl(\bolds{\beta}^{\star} \bigr)+ \mathbf{X}_{\mathcal{A}}'
\bolds{\Sigma} \bigl({\bolds{\beta}^{\star
}} \bigr)\mathbf{X}\bolds{
\Delta}+ \mathbf{R}_{\mathcal{A}}(\widetilde{\bolds{\Delta}}),
\]
where
$
\mathbf{R}_{\mathcal{A}}(\widetilde{\bolds{\Delta}})=\mathbf
{X}_{\mathcal{A}}'
(\bolds{\Sigma}(\bolds{\beta}^{\star}
+\widetilde{\bolds{\Delta}})-\bolds{\Sigma}({\bolds{\beta}^{\star
}})
)\mathbf{X}\bolds{\Delta}
$
with $\widetilde{\bolds{\Delta}}$ on the line segment joining $\mathbf{0}
$ and
$\bolds{\Delta}$.
Since $\bolds{\Delta}_{\mathcal{A}^c}=\mathbf{0}$, we have $\mathbf
{X}\bolds{\Delta}=\mathbf{X}
_{\mathcal{A}}\bolds{\Delta}_{\mathcal{A}}$.
By\vspace*{1.5pt} the mean-value theorem, we have\vspace*{1.5pt}
$
\|\mathbf{R}_{\mathcal{A}}(\widetilde{\bolds{\Delta}})\|_{\max}
\leq
\max_j\bolds{\Delta}_{\mathcal{A}}'\mathbf{X}_{\mathcal
{A}}'\operatorname{diag}
\{|\mathbf{x}_{(j)}|\circ|\bolds{\mu}''(\bar{\bolds{\beta}})|\}\mathbf
{X}_{\mathcal
{A}}\bolds{\Delta}_{\mathcal{A}}
$
for $\bar{\bolds{\beta}}$ being on the line segment joining
$\bolds{\beta}
^\star$ and $\bolds{\beta}^\star+ \widetilde{\bolds{\Delta}}$.
Using the fact that $|\psi'''(t)|=\theta(t)(1-\theta(t))|2\theta
(t)-1|\le\frac{1}4$ with $\theta(t)=(1+\exp(t))^{-1}$, we have
%
\begin{equation}
\bigl\|\mathbf{R}_{\mathcal{A}}(\widetilde{\bolds{\Delta}})\bigr\|_{\max} \leq
\frac
{n}4Q_1\cdot\|\bolds{\Delta}_{\mathcal{A}}
\|^2_{\ell_2} \leq\frac{n}4Q_1 s
r^2. \label{eq18}
\end{equation}

Notice that
\begin{eqnarray*}
F_{\mathcal{A}}(\bolds{\Delta}_{\mathcal{A}})&=& \bigl(\mathbf
{X}_{\mathcal{A}}' \bolds{\Sigma} \bigl(\bolds{
\beta}^{\star
} \bigr)\mathbf{X}_{\mathcal{A}} \bigr)^{-1} \bigl(
\mathbf{X}_{\mathcal{A}}'\mathbf{y}- \mathbf{X}_{\mathcal
{A}}'
\bolds{\mu} \bigl(\bolds{ \beta}^{\star
}+\bolds{\Delta} \bigr) \bigr)+
\bolds{\Delta}_{\mathcal{A}}
\\
&=& \bigl(\mathbf{X}_{\mathcal{A}}'\bolds{\Sigma} \bigl(\bolds{
\beta}^{\star
} \bigr)\mathbf{X}_{\mathcal
{A}} \bigr)^{-1}\cdot
\bigl(\mathbf{X}_{\mathcal{A}}'\mathbf{y}-\mathbf{X}_{\mathcal
{A}}'
\bolds{\mu} \bigl(\bolds{\beta}^{\star} \bigr)- \mathbf{R}_{\mathcal{A}}(
\widetilde{\bolds{\Delta}}) \bigr),
\end{eqnarray*}
we then use the triangle inequality to obtain
\begin{eqnarray*}
\bigl\|F_{\mathcal{A}}(\bolds{\Delta}_{\mathcal{A}})\bigr\|_{\max} &=& \bigl\| \bigl(
\mathbf{X}_{\mathcal{A}}'\bolds{\Sigma}(\bolds{\beta})\mathbf
{X}_{\mathcal{A}} \bigr)^{-1}\cdot \bigl(\mathbf{X}_{\mathcal{A}}'
\mathbf{y}- \mathbf{X}_{\mathcal{A}}'\bolds{\mu} \bigl(\bolds{
\beta}^{\star} \bigr)- \mathbf{R}_{\mathcal
{A}}(\widetilde{\bolds{
\Delta}}) \bigr)\bigr\|_{\max}
\\
&\le& Q_2\cdot \biggl(\biggl\|\frac{1}n\mathbf{X}_{\mathcal{A}}'
\bigl(\bolds{\mu} \bigl(\bolds{\beta}^{\star
} \bigr)-\mathbf{y} \bigr)
\biggr\|_{\max}+\frac{1}n\bigl\|\mathbf{R}_{\mathcal{A}}(\widetilde{\bolds{
\Delta}})\bigr\|_{\max
} \biggr).
\end{eqnarray*}
By (\ref{eq18}) and the definition of $r$, we have
$
\|F_{\mathcal{A}}(\bolds{\Delta}_{\mathcal{A}})\|_{\max} \leq
\frac{r}2+\frac{1}4Q_1Q_2sr^2 \leq r$.
This establishes the desired contraction (\ref{eq15}).

Next, we prove the upper bound for $\delta_1$. Recall that $\widehat
{\bolds{\Delta}}=\hat{\bolds{\beta}}{}^{\mathrm{oracle}}-\bolds{\beta
}^{\star}$.
Recall that
$
\ell_n(\bolds{\beta}) = \frac{1}{n}\sum_{i=1}^{n} \{
-y_i\mathbf{x}_i'\bolds{\beta}+\psi(\mathbf{x}_i'\bolds{\beta
}) \}$.
By a Taylor expansion,
%
\begin{equation}
\label{eq19} \qquad\nabla\ell_n \bigl(\hat{\bolds{\beta}}{}^{\mathrm{oracle}}
\bigr) = \nabla\ell_n \bigl(\bolds{\beta}^{\star} \bigr) +
\nabla^2\ell_n \bigl(\bolds{\beta}^{\star} \bigr)
\cdot\widehat{\bolds{\Delta}} + \bigl(\nabla^2\ell_n(
\tilde{\bolds{\beta}})-\nabla^2\ell_n \bigl(\bolds{\beta}
^{\star} \bigr) \bigr)\cdot\widehat{\bolds{\Delta}},
\end{equation}
where $\tilde{\bolds{\beta}}$ is on the line segment joining
$\hat{\bolds{\beta}}{}^{\mathrm{oracle}}$ and $\bolds{\beta}^{\star}$. Observe that
the first and
second derivatives of $\ell_n(\bolds{\beta})$ can be explicitly
written as
%
\begin{equation}
\label{eq20} \nabla\ell_n(\bolds{\beta})=\frac{1}n
\mathbf{X}' \bigl(\bolds{\mu}(\bolds{\beta})-\mathbf{y} \bigr)\quad
\mbox{and}\quad\nabla^2\ell_n(\bolds{\beta})=
\frac{1}n\mathbf{X}'\bolds{\Sigma}(\bolds{\beta})
\mathbf{X}.
\end{equation}
We define\vspace*{2pt}
$
\mathbf{R}(\bolds{\Delta})
= (\nabla^2\ell_n(\tilde{\bolds{\beta}})-\nabla^2\ell
_n(\bolds{\beta}
^{\star}) )\cdot\widehat{\bolds{\Delta}}
= \mathbf{X}' (\bolds{\Sigma}(\bolds{\beta}^{\star}+\bolds
{\Delta
})-\bolds{\Sigma}({\bolds{\beta}
^{\star}}) )\mathbf{X}\widehat{\bolds{\Delta}}$.
We rewrite\vspace*{1pt} $\mathbf{R}(\bolds{\Delta})$ as $(\mathbf{R}'_{\mathcal
{A}}(\bolds{\Delta}),\mathbf{R}
'_{\mathcal{A}^c}(\bolds{\Delta}))'$. Let $\widetilde{\bolds{\Delta
}}=\tilde{\bolds{\beta}}-\bolds{\beta}^{\star}$. Then, using $\widehat
{\bolds{\Delta}}_{\mathcal
{A}^c}=\mathbf{0}$, we have $\mathbf{X}\widehat{\bolds{\Delta}}=
\mathbf
{X}_{\mathcal{A}}
\widehat{\bolds{\Delta}}_{\mathcal{A}}$. Substituting this into (\ref
{eq19}), we obtain
%
\begin{equation}
\label{eq21} \nabla_{\mathcal{A}}\ell_n \bigl(\hat{\bolds{
\beta}}{}^{\mathrm{oracle}} \bigr)= \nabla_{\mathcal{A}}\ell_n \bigl(
\bolds{ \beta}^{\star} \bigr) +\frac{1}n\mathbf{X}_{\mathcal{A}}'
\bolds{\Sigma} \bigl(\bolds{\beta}^{\star
} \bigr)\mathbf{X}_{\mathcal
{A}}
\widehat{\bolds{\Delta}}_{\mathcal{A}} + \frac{1}{n} \mathbf
{R}_{\mathcal{A}}( \widetilde{\bolds{\Delta}})
\end{equation}
and
%
\begin{equation}
\label{eq22} \nabla_{\mathcal{A}^c}\ell_n \bigl(\hat{\bolds{
\beta}}{}^{\mathrm{oracle}} \bigr)= \nabla_{\mathcal{A}^c}\ell_n \bigl(
\bolds{ \beta}^{\star} \bigr)+\frac{1}n \mathbf{X}_{\mathcal{A}^c}'
\bolds{\Sigma} \bigl(\bolds{\beta}^{\star
} \bigr)\mathbf{X}_{\mathcal
{A}}
\widehat{\bolds{\Delta}}_{\mathcal{A}} + \frac{1}{n} \mathbf
{R}_{\mathcal{A}^c}( \widetilde{\bolds{\Delta}}).
\end{equation}
Using (\ref{eq20}) for $\bolds{\beta}^{\star}$ and $\nabla
_{\mathcal
{A}}\ell_n(\hat{\bolds{\beta}}{}^{\mathrm{oracle}})=\mathbf{0}$, we solve for
$\widehat{\bolds{\Delta}}_{\mathcal{A}}$ from (\ref{eq21}) and
substitute it into
(\ref{eq22})
to obtain
\begin{eqnarray*}
&&\nabla_{\mathcal{A}^c}\ell_n \bigl(\hat{\bolds{
\beta}}{}^{\mathrm{oracle}} \bigr)
\\
&&\qquad = \mathbf{X}_{\mathcal{A}^c}'\bolds{\Sigma} \bigl(
\bolds{\beta}^{\star
} \bigr)\mathbf{X}_{\mathcal{A}} \bigl(\mathbf{X}
_{\mathcal{A}}' \bolds{\Sigma} \bigl(\bolds{\beta}^{\star}
\bigr)\mathbf{X}_{\mathcal{A}} \bigr)^{-1} \biggl(-\frac{1}n
\mathbf{X}_{\mathcal{A}}' \bigl(\bolds{\mu} \bigl(\bolds{
\beta}^{\star
} \bigr)-\mathbf{y} \bigr)- \frac{1}{n}
\mathbf{R}_{\mathcal{A}}(\widetilde{\bolds{\Delta}}) \biggr)
\\
&&\quad\qquad{}+\frac{1}n\mathbf{X}_{\mathcal{A}^c}' \bigl(\bolds{\mu}
\bigl(\bolds{\beta}^{\star} \bigr)-\mathbf{y} \bigr)+ \frac{1}{n}
\mathbf{R}_{\mathcal{A}^c}(\widetilde{\bolds{\Delta}}).
\end{eqnarray*}

Recall that we have proved that (\ref{eq17}) holds under the condition
(\ref{eq16}). Now under the condition (\ref{eq16}) and the additional event
\[
\biggl\{ \bigl\|\nabla_{\mathcal{A}^c}\ell_n \bigl(\bolds{
\beta}^{\star} \bigr)\bigr\|_{\max
}<\frac
{a_1\lambda}2 \biggr\} \cap
\biggl\{ \bigl\|\nabla_{\mathcal{A}}\ell_n \bigl(\bolds{
\beta}^{\star} \bigr)\bigr\|_{\max}\le\frac{a_1\lambda}{4Q_3+2} \biggr\},
\]
we can follow the same lines of proof as in (\ref{eq18}) to show that
\[
\bigl\|\mathbf{R}(\widetilde{\bolds{\Delta}})\bigr\|_{\max} \le
\frac{n}4Q_1\|\widehat{\bolds{\Delta}}_{\mathcal{A}}
\|^2_{\ell_2} \le\frac{n}4Q_1sr^2,
\]
where $r=2Q_2\cdot\|\nabla_{\mathcal{A}}\ell_n(\bolds{\beta
}^{\star})\|
_{\max}$. Noticing that under condition (\ref{eq16})
\[
\frac{n}4Q_1sr^2= snQ_1Q_2^2
\cdot\bigl\|\nabla_{\mathcal{A}}\ell_n \bigl(\bolds{\beta}^{\star}
\bigr)\bigr\|_{\max}^2 \le n\cdot\bigl\|\nabla_{\mathcal{A}}
\ell_n \bigl(\bolds{\beta}^{\star} \bigr)\bigr\| _{\max},
\]
under the same event we have
\begin{eqnarray*}
\bigl\|\nabla_{\mathcal{A}^c}\ell_n \bigl(\hat{\bolds{
\beta}}{}^{\mathrm{oracle}} \bigr)\bigr\| &\le& Q_3\cdot \biggl(\bigl\|
\nabla_{\mathcal{A}} \ell_n \bigl(\bolds{\beta}^{\star
} \bigr)
\bigr\| _{\max}+ \frac{1}n\bigl\|\mathbf{R}_{\mathcal{A}}(\widetilde{
\bolds{\Delta}}) \bigr\|_{\max
} \biggr)
\\
&&{} +\bigl\|\nabla_{\mathcal{A}^c}\ell_n \bigl(\bolds{\beta}^{\star}
\bigr)\bigr\|_{\max
}+\frac
{1}n\bigl\|\mathbf{R}_{\mathcal{A}^c}(\widetilde{
\bolds{\Delta}})\bigr\|_{\max
}
\\
&\le& (2Q_3+1)\cdot\bigl\|\nabla_{\mathcal{A}}\ell_n \bigl(
\bolds{\beta}^{\star
} \bigr)\bigr\|_{\max} +\bigl\|\nabla_{\mathcal{A}^c}
\ell_n \bigl(\bolds{\beta}^{\star} \bigr)\bigr\|_{\max}
\\
&<& a_1\lambda.
\end{eqnarray*}
The desired probability bound can be obtained by using Proposition 4(a)
of \citet{fan2011}. This completes the proof of Theorem \ref{logistic}.

\subsection{Proof of Theorem \texorpdfstring{\protect\ref{logisticlasso}}{5}}\label{sec5.5}

The proof is relegated to a supplementary file [\citet{FZZ14supp}] 
for the sake of space constraint.

\subsection{Proof of Theorem \texorpdfstring{\protect\ref{glasso}}{6}}\label{sec5.6}

We\vspace*{1pt} first derive bound $\delta_2=\Pr(\|\widehat{\bolds{\Theta}}{}^{\mathrm{oracle}}_{\mathcal{A}}\|_{\min}\le a\lambda)$. A translation of
(\ref{eq3}) into the precision matrix estimation setting becomes
$
\widehat{\bolds{\Sigma}}{}^{\mathrm{oracle}}_\mathcal{A}=\widehat{\bolds{\Sigma
}}^n_{\mathcal{A}}$.
Let $\bolds{\Sigma}^{\bolds{\Delta}}=(\bolds{\Theta}^{\star
}+\bolds
{\Delta})^{-1}$ and
$r=2K_2\|\widehat{\bolds{\Sigma}}^n_{\mathcal{A}}-\bolds{\Sigma
}^{\star}_{\mathcal
{A}}\|_{\max}$. We define a map $F\dvtx \mathbb{B}(r)\subset\mathbb
{R}^{p^2}\rightarrow\mathbb{R}^{p^2}$ satisfying
$
F(\operatorname{vec}(\bolds{\Delta}))= ((F_{\mathcal{A}}(\operatorname{vec}(\bolds{\Delta
}_{\mathcal
{A}})))',\mathbf{0}' )'
$
with
%
\begin{equation}
\label{eq29} F_{\mathcal{A}} \bigl(\operatorname{vec}(\bolds{\Delta}_{\mathcal{A}}) \bigr)=
\bigl(\mathbf{H}^{\star}_{\mathcal
{A}\mathcal{A}} \bigr)^{-1}\cdot
\bigl(\operatorname{vec} \bigl(\bolds{\Sigma}^{\bolds{\Delta
}}_{\mathcal
{A}} \bigr)-\operatorname{vec} \bigl(
\widehat{\bolds{ \Sigma}}^n_{\mathcal{A}} \bigr) \bigr)+\operatorname{vec}(\bolds{
\Delta }_{\mathcal{A}})
\end{equation}
and
$
\mathbb{B}(r)=\{\bolds{\Delta}\dvtx \|\bolds{\Delta}_{\mathcal{A}}\|
_{\max}\le
r,\bolds{\Delta}_{\mathcal{A}^c}=\mathbf{0}\}$.
We will show that
%
\begin{equation}
\label{eq30} F \bigl(\mathbb{B}(r) \bigr)\subset\mathbb{B}(r)
\end{equation}
under the condition
%
\begin{equation}
\label{eq31} \bigl\|\widehat{\bolds{\Sigma}}^n_{\mathcal{A}}-\bolds{
\Sigma}^{\star
}_{\mathcal{A}}\bigr\| _{\max}< \min \biggl\{
\frac{1}{6K_1K_2d},\frac{1}{6K_1^3K_2^2d} \biggr\}.
\end{equation}
If\vspace*{2pt} (\ref{eq30}) holds, an application of the Brouwer's fixed-point
theorem yields a fixed-point $\widehat{\bolds{\Delta}}$ in\vspace*{1pt} the convex compact
set $\mathbb{B}(r)$ satisfying
$
\widehat{\bolds{\Delta}}_{\mathcal{A}^c}=\mathbf{0}
$
and
$
F_{\mathcal{A}}(\operatorname{vec}(\widehat{\bolds{\Delta}}_{\mathcal
{A}}))=\operatorname{vec}(\widehat{\bolds{\Delta}}_{\mathcal{A}})$.
Thus, $\widehat{\bolds{\Delta}}_{\mathcal{A}}=\widehat{\bolds{\Theta}}{}^{\mathrm{oracle}}_{\mathcal{A}}-\bolds{\Theta}^{\star}_{\mathcal{A}}$ by the
uniqueness and
%
\begin{equation}
\label{eq32} \bigl\|\widehat{\bolds{\Theta}}{}^{\mathrm{oracle}}-\bolds{
\Theta}^{\star}\bigr\|_{\max} = \|\widehat{\bolds{\Delta}}
\|_{\max} \le r.
\end{equation}

We now establish (\ref{eq30}). For any $\bolds{\Delta}\in\mathbb{B}(r)$,
by using (\ref{eq31}) we have
\[
\bigl\|\bolds{\Sigma}^{\star}\bolds{\Delta}\bigr\|_{\ell_\infty}\le
K_1\cdot\| \bolds{\Delta}\| _{\ell_1}\le K_1\cdot
\,dr=2K_1K_2d\cdot\bigl\|\widehat{\bolds{\Sigma}}^n_{\mathcal
{A}}-
\bolds{\Sigma}^{\star}_{\mathcal{A}}\bigr\|_{\max}<\tfrac{1}3.
\]
Thus,
$
\mathbf{J}=\sum_{j=0}^{\infty}(-1)^j(\bolds{\Sigma}^{\star}\bolds
{\Delta})^j
$
is a convergent matrix series of $\bolds{\Delta}$. Hence,
%
\begin{equation}
\label{{eq33}} \bolds{\Sigma}^{\bolds{\Delta}} = \bigl(\mathbf{I}+\bolds{
\Sigma}^{\star}\bolds{\Delta} \bigr)^{-1}\cdot\bolds{\Sigma}
^{\star} = \bolds{\Sigma}^{\star}-\bolds{\Sigma}^{\star}
\bolds{\Delta}\bolds{\Sigma}^{\star}+ \mathbf{R}^{\bolds{\Delta}},
\end{equation}
where $\mathbf{R}^{\bolds{\Delta}}=(\bolds{\Sigma}^{\star}\bolds
{\Delta
})^2\cdot\mathbf{J}\bolds{\Sigma}
^{\star}$.
Then it immediately yields that
%
\begin{eqnarray}\label{eq34}
&& \operatorname{vec} \bigl(\bolds{\Sigma}^{\bolds{\Delta}}_{\mathcal{A}} \bigr)-\operatorname{vec}
\bigl(\widehat{\bolds{\Sigma}} ^n_{\mathcal{A}} \bigr)
\nonumber\\[-8pt]\\[-8pt]
&&\qquad = \bigl(\operatorname{vec}
\bigl(\bolds{\Sigma}^{\star}_{\mathcal{A}} \bigr)-\operatorname{vec} \bigl(\widehat {
\bolds{ \Sigma}} ^n_{\mathcal{A}} \bigr) \bigr)- \operatorname{vec} \bigl(\bolds{
\Sigma}^{\star}\bolds{\Delta}\bolds{\Sigma} ^{\star
} \bigr)+\operatorname{vec}
\bigl(\mathbf{R}^{\bolds{\Delta}}_{\mathcal{A}} \bigr).\nonumber
\end{eqnarray}
Note that
$
\bolds{\Sigma}^{\star}\bolds{\Delta}\bolds{\Sigma}^{\star
}=(\bolds{\Sigma}^{\star
}\otimes
\bolds{\Sigma}^{\star})\cdot \operatorname{vec}(\bolds{\Delta})=\mathbf
{H}^{\star}\cdot
\operatorname{vec}(\bolds{\Delta})$,
and hence
\[
\operatorname{vec} \bigl(\bolds{\Sigma}^{\star}\bolds{\Delta}\bolds{
\Sigma}^{\star
}_{\mathcal{A}} \bigr) = \mathbf{H} ^\star_{\mathcal{A}\mathcal{A}}
\cdot \operatorname{vec}(\bolds{ \Delta}_{\mathcal{A}}).
\]
Now we follow the proof of Lemma 5 in \citet{ravikumar2008} to obtain
%
\begin{equation}
\label{eq35} \bigl\|\mathbf{R}^{\bolds{\Delta}}\bigr\|_{\max} = \max
_{(i,j)}\bigl|\mathbf{e}_i' \bigl( \bigl(
\bolds{\Sigma}^{\star}\bolds{\Delta} \bigr)^2\cdot\mathbf{J}
\bolds{\Sigma}^{\star} \bigr)\mathbf{e}_j\bigr| \le
\frac{3}2K^3_1\cdot d\|\Delta
\|_{\max}^2.
\end{equation}
Hence, a combination of (\ref{eq29}), (\ref{eq34}) and (\ref{eq35})
yields the contraction (\ref{eq30}), that is,
\begin{eqnarray*}
&&\bigl\|F_{\mathcal{A}} \bigl(\operatorname{vec}(\bolds{\Delta}_{\mathcal{A}}) \bigr)
\bigr\|_{\max}
\\
&&\qquad = \bigl\| \bigl(\mathbf{H}^{\star}_{\mathcal{A}\mathcal{A}} \bigr)^{-1}\cdot
\bigl( \bigl(\operatorname{vec} \bigl(\bolds{\Sigma}^{\star}_{\mathcal{A}} \bigr)-\operatorname{vec}
\bigl( \widehat{\bolds{\Sigma}} ^n_{\mathcal{A}} \bigr) \bigr)+\operatorname{vec}
\bigl( \mathbf{R}^{\bolds{\Delta}}_{\mathcal
{A}} \bigr) \bigr)\bigr\|_{\max
}
\\
&&\qquad \le K_2\cdot \bigl(\bigl\|\widehat{\bolds{\Sigma}}^n_{\mathcal{A}}-
\bolds{\Sigma}^{\star
}_{\mathcal{A}}\bigr\|_{\max}+\bigl\|
\mathbf{R}^{\bolds{\Delta}}\bigr\|_{\max} \bigr)
\\
&&\qquad \le r.
\end{eqnarray*}

Under the additional condition,
\[
\bigl\|\widehat{\bolds{\Sigma}}^n_{\mathcal{A}}-\bolds{
\Sigma}^{\star
}_{\mathcal{A}}\bigr\| _{\max}< \frac{1}{2K_2} \bigl(
\bigl\|\bolds{\Theta}^{\star}_{\mathcal{A}}\bigr\| _{\min
}-a\lambda \bigr)
\]
by (\ref{eq32}) and the definition or $r$, we have that
\begin{eqnarray*}
\bigl\|\widehat{\bolds{\Theta}}{}^{\mathrm{oracle}}_\mathcal{A}\bigr\|_{\min} &
\ge& \bigl\|\bolds{\Theta}^{\star}_{\mathcal{A}}\bigr\|_{\min}-\bigl\|\widehat{
\bolds{\Theta}}{}^{\mathrm{oracle}}-\bolds{\Theta}^{\star}\bigr\|_{\max}
\\
&=& \bigl\|\bolds{\Theta}^{\star}_{\mathcal{A}}\bigr\|_{\min}-2K_2
\cdot\bigl\| \widehat{\bolds{\Sigma}}^n_{\mathcal{A}}-\bolds{
\Sigma}^{\star}_{\mathcal
{A}}\bigr\|_{\max}
\\
&>&a\lambda.
\end{eqnarray*}
Thus,
\[
\delta_2 \le\Pr \biggl(\bigl\|\widehat{\bolds{\Sigma}}^n_{\mathcal{A}}-
\bolds{\Sigma}^{\star
}_{\mathcal{A}}\bigr\|_{\max}>\frac{1}{2K_2}
\min \biggl\{\frac{1}{3K_1d},\frac
{1}{3K_1^3K_2d}, \bigl\|\bolds{\Theta}^{\star}_{\mathcal{A}}
\bigr\|_{\min
}-a\lambda \biggr\} \biggr).
\]
An application of (\ref{eq8}) yields the bound on $\delta_2$.

We now bound $\delta_1$. Note that
$
\nabla_{\mathcal{A}^c}\ell_n(\widehat{\bolds{\Theta}}{}^{\mathrm{oracle}})
=
\widehat{\bolds{\Sigma}}^n_{\mathcal{A}^c}-\widehat{\bolds{\Sigma}}{}^{\mathrm{oracle}}_{\mathcal{A}^c}
$, and hence
%
\begin{equation}
\label{eq36} \bigl\|\nabla_{\mathcal{A}^c}\ell_n \bigl(\widehat{\bolds{
\Theta}}{}^{\mathrm{oracle}} \bigr)\bigr\| _{\max} \le\bigl\|\widehat{\bolds{
\Sigma}}^n_{\mathcal{A}^c}-\bolds{\Sigma}^{\star
}_{\mathcal
{A}^c}
\bigr\|_{\max}+ \bigl\|\widehat{\bolds{\Sigma}}{}^{\mathrm{oracle}}_{\mathcal{A}^c}-
\bolds{\Sigma}^{\star
}_{\mathcal{A}^c}\bigr\|_{\max}.
\end{equation}
Note\vspace*{2pt} $\|\widehat{\bolds{\Sigma}}^n_{\mathcal{A}^c}-\bolds{\Sigma
}^{\star
}_{\mathcal{A}^c}\|_{\max}$ is bounded by using (\ref{eq8}). Then we
only need to bound $\|\widehat{\bolds{\Sigma}}{}^{\mathrm{oracle}}_{\mathcal
{A}^c}-\bolds{\Sigma}^{\star}_{\mathcal{A}^c}\|_{\max}$. Recall
$\widehat{\bolds{\Delta}}=\widehat{\bolds{\Theta}}{}^{\mathrm{oracle}}-\bolds
{\Theta
}^{\star}$.
By (\ref{{eq33}}), we have
%
\begin{equation}
\label{eq37} \widehat{\bolds{\Sigma}}{}^{\mathrm{oracle}} = \bigl(\bolds{
\Theta}^{\star}+\widehat{\bolds{\Delta}} \bigr)^{-1} = \bigl(
\mathbf{I}+\bolds{\Sigma}^\star\bolds{\Delta} \bigr)^{-1}\cdot
\bolds{\Sigma}^\star= \bolds{\Sigma}^{\star}-\bolds{
\Sigma}^{\star}\widehat{\bolds{\Delta}}\bolds{\Sigma} ^{\star
}+
\widehat{\mathbf{R}},
\end{equation}
where $\widehat{\mathbf{R}}=(\bolds{\Sigma}^{\star}\widehat{\bolds
{\Delta}})^2\cdot
\widehat{\mathbf{J}}\bolds{\Sigma}^{\star}$ and $\widehat{\mathbf{J}}$
is defined
similar to $\mathbf{J}$ with $\bolds{\Delta}$ replaced by
$\widehat{\bolds{\Delta}}$.
Then $\widehat{\mathbf{J}}$ is a convergent matrix series under the condition
(\ref{eq31}). In terms of $\mathcal{A}$, we can equivalently write
(\ref{eq37}) as
\begin{eqnarray*}
\operatorname{vec} \bigl(\widehat{\bolds{\Sigma}}{}^{\mathrm{oracle}}_{\mathcal{A}} \bigr)-\operatorname{vec}
\bigl( \bolds{\Sigma}^{\star
}_{\mathcal{A}} \bigr) &=& -
\mathbf{H}^{\star}_{\mathcal{A}\mathcal{A}} \cdot \operatorname{vec}(\widehat{\bolds{
\Delta}}_{\mathcal{A}})+\operatorname{vec}(\widehat{ \mathbf{R}}_{\mathcal{A}}),
\\
\operatorname{vec} \bigl(\widehat{\bolds{\Sigma}}{}^{\mathrm{oracle}}_{\mathcal{A}^c} \bigr)-\operatorname{vec}
\bigl( \bolds{\Sigma}^{\star
}_{\mathcal{A}^c} \bigr) &=& -
\mathbf{H}^{\star}_{\mathcal{A}^c\mathcal{A}} \cdot \operatorname{vec}(\widehat{\bolds{\Delta}}
_{\mathcal{A}})+\operatorname{vec}(\widehat{ \mathbf{R}}_{\mathcal{A}^c}),
\end{eqnarray*}
where we use the fact that $\widehat{\bolds{\Delta}}_{\mathcal
{A}^c}=\mathbf{0}
$. Solving $\operatorname{vec}(\widehat{\bolds{\Delta}}_{\mathcal{A}})$ from the first
equation and substituting it into the second equation, we obtain
\begin{eqnarray*}
&&\operatorname{vec} \bigl(\widehat{\bolds{\Sigma}}{}^{\mathrm{oracle}}_{\mathcal{A}^c} \bigr)-\operatorname{vec}
\bigl(\bolds{\Sigma}^{\star
}_{\mathcal{A}^c} \bigr)
\\
&&\qquad = \mathbf{H}^{\star}_{\mathcal{A}^c\mathcal{A}} \bigl(\mathbf{H}^{\star
}_{\mathcal
{A}\mathcal{A}}
\bigr)^{-1}\cdot \bigl(\operatorname{vec} \bigl(\widehat{\bolds{\Sigma}}{}^{\mathrm{oracle}}_{\mathcal{A}} \bigr)-\operatorname{vec} \bigl(\bolds{\Sigma}^{\star}_{\mathcal
{A}}
\bigr)-\operatorname{vec}(\widehat{\mathbf{R}}_{\mathcal{A}}) \bigr)+\operatorname{vec}(\widehat {\mathbf{R}}
_{\mathcal{A}^c}).
\end{eqnarray*}

Recall (\ref{eq35}) holds under condition (\ref{eq31}). Thus, we have
\[
\|\widehat{\mathbf{R}}\|_{\max} \le\tfrac{3}2K^3_1
\cdot d\|\widehat\Delta\|_{\max}^2 = 6K^3_1K^2_2
d\cdot\bigl\|\widehat{\bolds{\Sigma}}_{\mathcal{A}}^n-\bolds{\Sigma}
_{\mathcal{A}}{}^{\star}\bigr\|_{\max} \le\bigl\|\widehat{\bolds{
\Sigma}}_{\mathcal{A}}^n-\bolds{\Sigma}_{\mathcal
{A}}{}^{\star}
\bigr\| _{\max}.
\]
Under the extra event
$
\{\|\widehat{\bolds{\Sigma}}_{\mathcal{A}^c}^n-\bolds{\Sigma
}_{\mathcal
{A}^c}^{\star}\|_{\max}<
\frac{a_1\lambda}{2} \}\cap
\{\|\widehat{\bolds{\Sigma}}^n_{\mathcal{A}}-\bolds{\Sigma
}^{\star
}_{\mathcal{A}}\|_{\max}\le
{\frac{a_1\lambda}{4K_3+2}} \}$,
we derive the desired upper bound for (\ref{eq36}) by using the
triangular inequality,
\begin{eqnarray*}
\bigl\|\nabla_{\mathcal{A}^c}\ell_n \bigl(\widehat{\bolds{
\Theta}}{}^{\mathrm{oracle}} \bigr)\bigr\| _{\max} &\le& \bigl\|\widehat{\bolds{
\Sigma}}^n_{\mathcal{A}^c}-\bolds{\Sigma}^{\star
}_{\mathcal
{A}^c}
\bigr\|_{\max}+\bigl\|\bolds{\Sigma}^{\star}_{\mathcal{A}^c}-\widehat{\bolds{
\Sigma}}{}^{\mathrm{oracle}}_{\mathcal{A}^c}\bigr\|_{\max}
\\
&\le& \frac{a_1\lambda}{2}+(2K_3+1)\cdot\bigl\|\widehat{\bolds{
\Sigma}}_{\mathcal
{A}}^n-\bolds{\Sigma}_{\mathcal{A}}{}^{\star}
\bigr\|_{\max}
\\
&<& a_1\lambda.
\end{eqnarray*}
Therefore,
\begin{eqnarray*}
\delta_1 &\leq& \Pr \biggl\{\bigl\|\widehat{\bolds{\Sigma}}^n_{\mathcal
{A}}-
\bolds{\Sigma} ^{\star}_{\mathcal{A}}\bigr\|_{\max} \ge\min \biggl\{
\frac{1}{6K_1K_2d},\frac
{1}{6K_1^3K_2^2d}, \frac{a_1\lambda}{4K_3+2} \biggr\} \biggr\}
\\
& &{} + \Pr \biggl\{\bigl\|\widehat{\bolds{\Sigma}}_{\mathcal{A}^c}^n-\bolds
{\Sigma}_{\mathcal
{A}^c}^{\star}\bigr\|_{\max} > \frac{a_1\lambda}{2}
\biggr\}.
\end{eqnarray*}
An application of (\ref{eq8}) yields $\delta_1^G$. This completes the
proof of Theorem \ref{glasso}.

\subsection{Proof of Theorem \texorpdfstring{\protect\ref{quantile}}{7}}\label{sec5.7}
First, we bound $\delta_2=\Pr(\|\hat{\bolds{\beta}}_{\mathcal
{A}}{}^{\mathrm{oracle}}\|_{\min}\le a\lambda)$. To this end, we let
\[
\mathbb{B}(r)= \bigl\{\bolds{\Delta}\in\mathbb{R}^p\dvtx \|\bolds{
\Delta }_{\mathcal{A}}\| _{\ell_2}\le r,\bolds{\Delta}_{\mathcal{A}^c}=
\mathbf{0} \bigr\}
\]
with\vspace*{1pt} $r=\|\bolds{\beta}^\star_{\mathcal{A}}\|_{\min}-a\lambda$, and
$\partial\mathbb{B}(r) = \{\bolds{\Delta}\in\mathbb{R}^p\dvtx \|\bolds
{\Delta}
_{\mathcal{A}}\|_{\ell_2}= r,\bolds{\Delta}_{\mathcal{A}^c}=\mathbf{0}
\}$. We define
$
F(\bolds{\Delta})
=\ell_n(\bolds{\beta}^\star+\bolds{\Delta})-\ell_n(\bolds{\beta
}^\star)
$
and $\widehat{\bolds{\Delta}}=\arg\min_{\bolds{\Delta}\dvtx  \bolds
{\Delta}_{\mathcal
{A}^c}=\mathbf{0}} F(\bolds{\Delta})$. Then $\widehat{\bolds{\Delta
}}=\hat{\bolds{\beta}}{}^{\mathrm{oracle}}-\bolds{\beta}^\star$. Since\vspace*{1pt} $F(\widehat{\bolds{\Delta}})\le
F(\mathbf{0})=0$
holds by definition, the convexity of $F(\bolds{\Delta})$ yields that
$
\Pr(\|\widehat{\bolds{\Delta}}_{\mathcal{A}}\|_{\ell_2} \le r )
\ge
\Pr(\inf_{\bolds{\Delta}\in\partial\mathbb{B}(r)}F(\bolds
{\Delta})>0)$,
and\vspace*{2pt} thus
$
\delta_2
\le
1-\Pr(\|\widehat{\bolds{\Delta}}_{\mathcal{A}}\|_{\ell_2} \le r )
\le
1-\Pr(\inf_{\bolds{\Delta}\in\partial\mathbb{B}(r)}F(\bolds
{\Delta})>0)$.
In what follows, it suffices to bound $\Pr(\inf_{\bolds{\Delta}\in
\partial
\mathbb{B}(r)}F(\bolds{\Delta})>0)$.

By the definition of $\rho_{\tau}(\cdot)$ and $y_i=\mathbf
{x}_i'\bolds{\beta}^\star+\varepsilon_i$, we can rewrite $F(\bolds
{\Delta})$ as
\begin{eqnarray*}
F(\bolds{\Delta}) &=& \frac{1}{n}\sum_i
\bigl\{\rho_{\tau} \bigl(y_i-\mathbf{x}_i'
\bigl(\bolds{\beta} ^\star+\bolds{\Delta} \bigr) \bigr)-
\rho_{\tau} \bigl(y_i-\mathbf{x}_i'
\bolds{\beta}^\star \bigr) \bigr\}
\\
&=& \frac{1}{n}\sum_i \bigl\{
\rho_{\tau} \bigl(\varepsilon_i-\mathbf{x}_i'
\bolds{\Delta} \bigr)- \rho_{\tau}(\varepsilon_i) \bigr\}
\\
&=& \frac{1}{n}\sum_i \bigl\{
\mathbf{x}_i'\bolds{\Delta}\cdot(I_{\{\varepsilon_i\le0\}}-
\tau) + \bigl(\mathbf{x}_i'\bolds{\Delta}-
\varepsilon_i \bigr)\cdot(I_{\{\varepsilon
_i\le\mathbf{x}_i'\bolds{\Delta}\}}-I_{\{\varepsilon_i\le0\}}) \bigr
\}.
\end{eqnarray*}
Next, we bound $I_1= \frac{1}{n}\sum_i\mathbf{x}_i'\bolds{\Delta
}\cdot
(I_{\{\varepsilon_i\le0\}}-\tau)$ and $I_2=F(\bolds{\Delta})-I_1$,
respectively.

To bound $I_1$, we use the Cauchy--Schwarz inequality to obtain
\[
\bigl|\mathbf{x}_i'\bolds{\Delta}\cdot I_{\{\varepsilon_i\le0\}}\bigr|
\le\bigl|\mathbf{x}_i'\bolds{\Delta}\bigr| = |
\mathbf{x}_{i\mathcal{A}}\bolds{\Delta}_{\mathcal{A}}| \le\|
\mathbf{x}_{i\mathcal{A}}\|_{\ell_2} \cdot\|\bolds{\Delta} _{\mathcal{A}}
\|_{\ell_2} = M^{1/2}_{\mathcal{A}} s^{1/2}r.
\]
Since $E[\mathbf{x}_i'\bolds{\Delta}\cdot I_{\{\varepsilon_i\le0\}
}]=\mathbf{x}_i'\bolds{\Delta}\cdot\Pr(\varepsilon_i\le0) =
\mathbf{x}_i'\bolds{\Delta}\cdot\tau$, we can apply the Hoeffding's
inequality to bound $I_1$ as follows:
%
\begin{equation}
\label{eq38} \Pr \biggl(|I_1|>\frac{1}6
\lambda_{\min} f_{\min} r^2 \biggr) 
\le2\exp \biggl(-\frac{\lambda_{\min}^2 f_{\min}^2}{72M_{\mathcal
{A}}\cdot
s}\cdot nr^2
\biggr).
\end{equation}

Now we bound $I_2$. Using Knight's identity [\citet{knight1998}],
we write $I_2$ as
\[
I_2 = \frac{1}n\sum_i \bigl(
\mathbf{x}_i'\bolds{\Delta}-\varepsilon_i
\bigr)\cdot(I_{\{
\varepsilon_i\le\mathbf{x}_i'\bolds{\Delta}\}}-I_{\{\varepsilon
_i\le
0\}}) = \frac{1}n\sum
_i\int_0^{\mathbf{x}_i'\bolds{\Delta}}
(I_{\{
\varepsilon
_i\le s\}}-I_{\{\varepsilon_i\le0\}}) \,ds.
\]
Note that each term in the summation of $I_2$ is uniformly bounded,
that is,
\[
\biggl|\int_0^{\mathbf{x}_i'\bolds{\Delta}} (I_{\{\varepsilon_i\le s\}
}-I_{\{
\varepsilon_i\le0\}})
\,ds\biggr| \le\biggl|\int_0^{\mathbf{x}_i'\bolds{\Delta}} 1\,ds\biggr| \le\bigl|
\mathbf{x}_i'\bolds{\Delta}\bigr| \le M_{\mathcal{A}}{}^{1/2}s^{1/2}
r.
\]
Then we can use the Hoeffding's inequality to bound $I_2-E[I_2]$ as
%
\begin{equation}
\label{eq39} \Pr \biggl(\bigl|I_2-E[I_2]\bigr|>
\frac{1}6 \lambda_{\min} f_{\min} r^2 \biggr)
\le2\exp \biggl(- \frac{\lambda_{\min}^2 f_{\min}^2}{72M_{\mathcal{A}}\cdot
s}\cdot nr^2 \biggr).
\end{equation}
Next, we apply Fubini's theorem and mean-value theorem to derive that
\begin{eqnarray*}
E[I_2] 
&=& \frac{1}n
\sum_i\int_0^{\mathbf{x}_i'\bolds{\Delta}}
\bigl(F_i(s)-F_i(0) \bigr) \,ds
\\
&=& \frac{1}n\sum_i\int
_0^{\mathbf{x}_i'\bolds{\Delta}} f_i \bigl(\xi(s) \bigr)\cdot
s \,ds,
\end{eqnarray*}
where $\xi(s)$ is on the line segment between $0$ and $s$. By the
assumption\vspace*{1.5pt} of Theorem~\ref{quantile}, $|\xi(s)|\le|\mathbf
{x}_i'\bolds{\Delta}|\le M_{\mathcal{A}}{}^{1/2}s r^{1/2}\le u_0$
holds,\vspace*{1pt} and
then, by condition \textup{(C3)}, we have $f_i(\xi(s))\ge f_{\min}$ for any
$i$. Using this fact, it is easy to obtain
\[
E[I_2] \ge\frac{1}n\sum_i
\int_0^{\mathbf{x}_i'\bolds{\Delta}} f_{\min
}\cdot s \,ds =
\frac{1}{2n}f_{\min}\sum_i (
\mathbf{x}_{i\mathcal{A}}\bolds{\Delta} _{\mathcal{A}})^2
\ge\frac{1}{2}
\lambda_{\min} f_{\min} r^2.
\]
This together with (\ref{eq38}) and (\ref{eq39}) proves that under
the event
\[
\bigl\{|I_1|\le\tfrac{1}6 \lambda_{\min}
f_{\min} r^2 \bigr\}\cup \bigl\{ \bigl|I_2-E[I_2]\bigr|
\le\tfrac{1}6 \lambda_{\min} f_{\min} r^2
\bigr\},
\]
we have
$
F(\bolds{\Delta}) =I_1+I_2
\ge- |I_1| + E[I_2] + (I_2-E[I_2])
\ge\frac{1}{6} \lambda_{\min} f_{\min} r^2 >0$.
Therefore, an application of the union bound yields the desired bound
$\delta^Q_2$.

In the sequel, we bound $\delta_1=\Pr(\Vert\nabla_{\mathcal
{A}^c}\ell_n(\hat{\bolds{\beta}}{}^{\mathrm{oracle}})\Vert_{\max} \ge
a_1\lambda
)$. A translation of~(\ref{eq3}) implies the subgradient optimality
condition $\nabla_{j}\ell_n(\hat{\bolds{\beta}}{}^{\mathrm{oracle}})=0$ for
$j\in
\mathcal{A}$. For ease of notation, let $\hat\varepsilon_i = y_i -
\mathbf{x}_{i}' \hat{\bolds{\beta}}{}^{\mathrm{oracle}}$ and $\mathcal{Z}=\{i\dvtx
\hat\varepsilon_i =0\}$. Now we can rewrite the subgradient $\nabla
_{j}\ell_n(\hat{\bolds{\beta}}{}^{\mathrm{oracle}})$ as
\begin{eqnarray*}
\nabla_{j}\ell_n \bigl(\hat{\bolds{\beta}}{}^{\mathrm{oracle}}
\bigr) &=& \frac{1}n\sum_{i\notin\mathcal{Z}}
x_{ij}\cdot(I_{\{\hat
\varepsilon_i\le0\}}-\tau) -\frac{1}n\sum
_{i\in\mathcal{Z}} x_{ij}\cdot\hat z_i
\\
&=& \frac{1}n\sum_{i} x_{ij}
\cdot(I_{\{\hat\varepsilon_i\le0\}}-\tau) -\frac{1}n\sum_{i\in\mathcal{Z}}
x_{ij}\cdot(\hat z_i +1-\tau)
\\
&=& I_{3j}+I_{4j},
\end{eqnarray*}
where $\hat z_i$ $\in[\tau-1,\tau]$ ($i\in\mathcal{Z}$) satisfies
the subgradient optimality condition. To bound $\nabla_{j}\ell
_n(\hat{\bolds{\beta}}{}^{\mathrm{oracle}})$ for $j\in\mathcal{A}^c$, a key
observation is that the quantile regression for $\hat{\bolds{\beta}}{}^{\mathrm{oracle}}$ exactly interpolates $s$ observations, that is, $|\mathcal
{Z}|=s$. Please see Section~2.2 of \citet{koenker2005} for more
details. Then it is easy to derive
\[
\max_{j\in\mathcal{A}^c}|I_{4j}| \le\max_{j\in\mathcal{A}^c}
\frac{1}n\sum_{i\in\mathcal{Z}} |x_{ij}| \cdot
\bigl(\max\{1-\tau, \tau\}+1-\tau \bigr) \le2m_{\mathcal{A}^c}\cdot
\frac{s}n \le\frac{a_1}4\lambda.
\]
Using this bound for $|I_{4j}|$, we can further bound $\delta_1$ as
%
\begin{equation}
\label{eq40} \delta_1 \le\Pr \Bigl(\max_{j\in\mathcal{A}^c}|I_{3j}+I_{4j}|
\ge a_1\lambda \Bigr) \le\Pr \biggl(\max_{j\in\mathcal{A}^c}|I_{3j}|
\ge\frac{3a_1}4\lambda \biggr). 
\end{equation}

Now we only need to bound $\max_{j\in\mathcal{A}^c}|I_{3j}|$. Note
that we rewrite $I_{3j}$ as
\begin{eqnarray*}
I_{3j} 
&=& \frac{1}n\sum
_{i} x_{ij} \bigl((I_{\{\hat\varepsilon_i\le0\}
}-I_{\{\varepsilon_i\le0\}})
- E[I_{\{\hat\varepsilon_i\le0\}
}-I_{\{\varepsilon_i\le0\}}] \bigr)
\\
&&{} + \frac{1}n\sum_{i} x_{ij}
E[I_{\{\hat\varepsilon_i\le0\}}-I_{\{
\varepsilon_i\le0\}}] + \frac{1}n\sum
_{i} x_{ij}(I_{\{\varepsilon
_i\le0\}}-\tau)
\\
&=& I_{3j.1}+I_{3j.2}+I_{3j.3}.
\end{eqnarray*}
Next, we define $\varepsilon_i(\mathbf{t}) = y_i - \mathbf
{x}_{i\mathcal
{A}}' \mathbf{t}$. Then $\hat\varepsilon_i=\varepsilon
_i(\hat{\bolds{\beta}}_{\mathcal{A}}{}^{\mathrm{oracle}})$ holds by
definition. We
also introduce $\mathbb{B}_{\star}=\{\mathbf{t}\in\mathbb{R}^s\dvtx \|
\mathbf{t}
-\bolds{\beta}_{\mathcal{A}}{}^{\star}\|_{\ell_2}\le r_1=6\sqrt {M_{\mathcal
{A}}s\log n/n}\}$. As long as $\hat{\bolds{\beta}}_{\mathcal
{A}}{}^{\mathrm{oracle}}\in\mathbb{B}_{\star}$ holds, due to the mean-value
theorem, we have
\begin{eqnarray*}
|I_{3j.2}| &\le& \sup_{\mathbf{t}\in\mathbb{B}_{\star}(r)}\biggl|\frac{1}n\sum
_{i} x_{ij} E[I_{\{\varepsilon_i(\mathbf{t})\le0\}}-I_{\{\varepsilon
_i\le0\}
}]\biggr|
\\
&\le& \sup_{\mathbf{t}\in\mathbb{B}_{\star}(r)}\frac{m_{\mathcal
{A}^c}}n\sum
_{i} \bigl|F_i \bigl(\mathbf{x}_{i\mathcal{A}}'
\bigl(\mathbf{t}-\bolds{\beta}^{\star
}_{\mathcal{A}} \bigr)
\bigr)-F_i(0)\bigr|
\\
&\le& \sup_{\mathbf{t}\in\mathbb{B}_{\star}(r)} \frac{m_{\mathcal
{A}^c}}n\sum
_{i}\bigl|f_i \bigl(\xi_i(\mathbf{t})
\bigr)\cdot\mathbf{x}_{i\mathcal{A}}' \bigl(\mathbf{t} -\bolds{
\beta}^{\star}_{\mathcal{A}} \bigr)\bigr|,
\end{eqnarray*}
where\vspace*{1pt} $\xi_i(\mathbf{t})$ is on the line segment between $0$ and
$\mathbf{x}_{i\mathcal{A}}'(\mathbf{t}-\bolds{\beta}^{\star
}_{\mathcal
{A}})$. Note that $\sup_{\mathbf{t}\in\mathbb{B}_{\star
}(r)}|\mathbf
{x}_{i\mathcal{A}}'(\mathbf{t}-\bolds{\beta}^{\star}_{\mathcal
{A}})|\le
\|
\mathbf{x}_{i\mathcal{A}}\|_{\ell_2}\cdot\|\mathbf{t}-\bolds
{\beta
}^{\star
}_{\mathcal{A}}\|_{\ell_2}\le M_{\mathcal{A}}{}^{1/2} s^{1/2}r_1\le
u_0$.\vspace*{1pt} Then $\sup_{\mathbf{t}\in\mathbb{B}_{\star}(r)}|f_i(\xi
_i(\mathbf{t}
))|\le f_{\max}$ holds by condition \textup{(C3)}. Thus, we have
%
\begin{equation}
\label{eq41} \max_{j\in\mathcal{A}^c}|I_{3j.2}| \le
\frac{m_{\mathcal{A}^c}}n \cdot n f_{\max} \cdot M_{\mathcal
{A}}{}^{1/2}
s^{1/2}r_1 \le\frac{a_1}{4}\lambda.
\end{equation}

Let\vspace*{1pt} $\gamma_i(\mathbf{t})=I_{\{\varepsilon_i(\mathbf{t})\le0\}
}-I_{\{\varepsilon
_i\le0\}} - E[I_{\{\varepsilon_i(\mathbf{t})\le0\}}-I_{\{
\varepsilon_i\le
0\}}]$ and $I_{3j.1}(\mathbf{t})=\break \frac{1}n\sum_{i} x_{ij}\gamma
_i(\mathbf{t})$.
Again\vspace*{1pt} if $\hat{\bolds{\beta}}_{\mathcal{A}}{}^{\mathrm{oracle}}\in\mathbb
{B}_{\star}$, we have $|I_{3j.1}|\le\sup_{\mathbf{t}\in\mathbb
{B}_{\star
}(r)}|I_{3j.1}(\mathbf{t})|$. \mbox{Together} with this, we combine (\ref{eq40}),
(\ref{eq41}) and the union bound to obtain
%
\begin{eqnarray}\label{eq42}
\delta_1 &\le& \Pr \bigl(\hat{\bolds{\beta}}_{\mathcal{A}}{}^{\mathrm{oracle}}
\notin\mathbb{B}_{\star} \bigr) + \Pr \biggl(\max_{j\in\mathcal{A}^c}
\sup_{\mathbf
{t}\in
\mathbb{B}_{\star}(r)}\bigl|I_{3j.1}(\mathbf{t})\bigr|> \frac{a_1\lambda
}4
\biggr)
\nonumber\\[-8pt]\\[-8pt]
&&{}+ \Pr \biggl(\max_{j\in\mathcal{A}^c}|I_{3j.3}|>
\frac{a_1\lambda
}4 \biggr).\nonumber
\end{eqnarray}
Note that $\hat{\bolds{\beta}}_{\mathcal{A}}{}^{\mathrm{oracle}}\in\mathbb
{B}_{\star}$ holds under the event $\{ \| \hat{\bolds{\beta}}{}^{\mathrm{oracle}}_{\mathcal{A}} - \bolds{\beta}^\star_{\mathcal{A}}\|
_{\ell_2}\le
r_1\}$. Then we can combine (\ref{eq38}) and (\ref{eq39}) to derive that
%
\begin{equation}
\label{eq43} \Pr \bigl(\hat{\bolds{\beta}}_{\mathcal{A}}{}^{\mathrm{oracle}}\notin
\mathbb{B}_{\star} \bigr) \le\Pr \bigl( \bigl\| \hat{\bolds{
\beta}}{}^{\mathrm{oracle}}_{\mathcal{A}} - \bolds{\beta}^\star
_{\mathcal{A}}\bigr\|_{\ell_2} > r_1 \bigr) 
= 4n^{-1/2}.
\end{equation}
By the assumption of $\lambda$, we use Lemma A.3 of \citet
{wang2012} to obtain
%
\begin{equation}
\Pr \biggl( \sup_{\mathbf{t}\in\mathbb{B}_{\star}(r)} \biggl|\frac{1}n\sum
_{i} x_{ij} \gamma_i(
\mathbf{t})\biggr| > \frac
{a_1}4\lambda \biggr) 
\le C_1 (p-s)\exp \biggl(-\frac{a_1n\lambda}{104m_{\mathcal{A}^c}}
\biggr),\label{eq44}
\end{equation}
where $C_1>0$ is a fixed constant that does not depend on $n$, $p$,
$s$, $m_{\mathcal{A}^c}$ and $M_{\mathcal{A}}$. Furthermore, we use
the Hoeffding's inequality to bound $I_{3j.3}$ as
%
\begin{equation}
\Pr \biggl(\max_{j\in\mathcal{A}^c}|I_{3j.3}|>\frac{a_1}4
\lambda \biggr) \le2(p-s)\cdot\exp \biggl(-\frac{a_1^2n\lambda
^2}{32m_{\mathcal{A}^c}^2} \biggr).
\label{eq45}
\end{equation}

Therefore, we can combine (\ref{eq42}), (\ref{eq43}), (\ref{eq44})
and (\ref{eq45}) to obtain the desired probability bound $\delta_1^Q$
for $\delta_1$. This complete the proof of Theorem \ref{quantile}.

\section*{Acknowledgements}
We thank the Editor, Associate Editor and referees for their helpful
comments that improved an earlier version of this paper.

\begin{supplement}
\stitle{Supplement to ``Strong oracle optimality of folded concave penalized estimation''}
\slink[doi]{10.1214/13-AOS1198SUPP} 
\sdatatype{.pdf}
\sfilename{AOS1198\_supp.pdf}
\sdescription{In this supplementary note, we give the complete proof
of Theorem \ref{logisticlasso} and some comments on the simulation studies.}
\end{supplement}


%

\printaddresses

\end{document}